\documentclass[twocolumn]{autart}

\usepackage{floatrow}
\usepackage{picins}
\usepackage{algorithmic}
\usepackage{hyperref}
\usepackage{graphicx}
\usepackage[cmex10]{amsmath}
\usepackage{array}
\usepackage{mdwtab}
\usepackage{graphics}
\usepackage{epsfig}
\usepackage{mathptmx}
\usepackage{times}
\usepackage{amsmath}
\usepackage{amssymb}
\usepackage{url}
\usepackage{lscape}
\usepackage{color}
\usepackage{epsfig}
\usepackage{IEEEtrantools}
\usepackage{caption}
\usepackage{cite}
\usepackage{balance}
\usepackage{dblfloatfix}
\usepackage{eucal}
\usepackage{nicefrac}

\begin{document}

\begin{frontmatter}

\title{Reference Governor-Based Fault-Tolerant Constrained Control\thanksref{footnoteinfo}}

\thanks[footnoteinfo]{This research has been supported by National Science Foundation under award numbers ECCS-1931738 and ECCS-1932530.}

\author[WASHUESE]{Mehdi Hosseinzadeh\corauthref{Corres}}\ead{mehdi.hosseinzadeh@ieee.org},
\author[UoM]{Ilya Kolmanovsky}\ead{ilya@umich.edu}, \author[WASHUCSE]{Sanjoy Baruah}\ead{baruah@wustl.edu}, \author[WASHUESE]{Bruno Sinopoli}\ead{bsinopoli@wustl.edu}

\corauth[Corres]{Corresponding author.}

\address[WASHUESE]{Department of Electrical and Systems Engineering, Washington University in St. Louis, St. Louis, MO 63130, USA}
\address[UoM]{Department of Aerospace Engineering, The University of Michigan, Ann Arbor, MI 48109, USA}
\address[WASHUCSE]{Department of Computer Science and  Engineering, Washington University in St. Louis, St. Louis, MO 63130, USA}

\begin{abstract}
This paper presents a fault-tolerant control scheme for constrained linear systems. First, a new variant of the Reference Governor (RG) called At Once Reference Governor (AORG) is introduced. The AORG is distinguished from the conventional RG by computing the Auxiliary Reference (AR) sequence so that to optimize performance over a prescribed time interval instead of only at the current time instant; this enables the integration of the AORG  with fault detection schemes. In particular, it is shown that, when the AORG is combined with a Multi-Model Adaptive Estimator (MMAE), the AR sequence can be determined such that the tracking properties are guaranteed and constraints are satisfied at all times, while the detection performance is optimized, i.e., faults can be detected with a high probability of correctness. In addition a reconfiguration scheme is presented that ensures system \textit{viability} despite the presence of faults based on recoverable sets. Simulations on a Boeing 747-100 aircraft model are carried out to evaluate the effectiveness of the AORG scheme in enforcing constraints and tracking the desired roll and side-slip angles. The effectiveness of the presented fault-tolerant control scheme in maintaining the airplane viability in the presence of damaged vertical stabilizer is also demonstrated. 
\end{abstract}


\end{frontmatter}

%%%%%%%%%%%%%%%%%%%%%%%%%%%%%%%%%%%%%%%%%%%%%%%%%%%%%%%%%%%%%%%%%%%%%%%%%%%%%%%%
\section{Introduction}
The satisfaction of constraints (e.g. operational limits and actuator range and rate limits) is a crucial requirement for the control of many real-world systems. There are two typical choices to ensure constraint satisfaction. One choice is to design the controller within the model predictive control framework \cite{Mayne2000,Domahidi2012}. The other choice is to decouple the problem of the stabilization of the system from the problem of satisfying the constraints \cite{Kalabic2020}. In particular, a prestabilized system can be augmented with an \textit{add-on} unit called Reference Governor (RG) that, whenever necessary, modifies the reference signal to ensure constraint satisfaction  \cite{Bemporad1998,Gilbert1999,Garone2016_1}. Notably, a novel scheme called Explicit RG  has been introduced recently \cite{Hosseinzadeh2018ERG,Hosseinzadeh2019ECC,Nicotra2018,Hosseinadeh2019,Hosseinzadeh2020}, which deals with constrained reference tracking without resorting to on-line optimization.

Equipment faults/failures are the main source of industrial safety hazards \cite{Chiang2001,Raimondo2013,HosseinzadehWind}. As a result, designing a suitable fault-tolerant control scheme to mitigate the impacts of faults on the stability and performance of the systems has gained a great attention in recent years, e.g., \cite{Ashari2012,Blanke2016,HosseinzadehSmartGrid}. The fault-tolerant control schemes presented in the literature typically consist of two units \cite{Li2005,Seron2012}: 1) a fault detection unit, which detects the presence of a fault and identifies its nature, and 2) a control reconfiguration strategy, which modifies the control law to continue operating the system with potentially decreased/degraded functionality/availability despite the presence of the fault.

One key issue that is overlooked in most of existing literature is the system viability, which is characterized by measures of operational capability and satisfaction of operating constraints. Note that when a fault occurs, in many real-world applications, the most immediate objective is not to recover asymptotic properties (e.g. stability), but to ensure that the constraints are not violated during the transient. Indeed, the violation of constraints may have catastrophic consequences, making it impossible to recover a safe operation.

The need for fault-tolerant constrained control has been recognized in \cite{Puncochar2015}, where a control scheme has been presented which ensures constraint satisfaction despite the presence of faults, while optimizing control and detection performances. The scheme presented in \cite{Puncochar2015} applies the control sequence in open loop, which may make the system vulnerable to disturbances or model mismatch. Some fault-tolerant constrained control schemes based on model predictive control \cite{Maciejowski1999,Camacho2010} and RG \cite{Garone2016_1} have been presented in the literature as well. {\color{black} In particular, \cite{Maciejowski2003} considers the application of an MPC-based fault-tolerant control to deal with failures in both engines of a Boeing 747-200F. In \cite{Riverso2014,Boem2020} a distributed MPC fault-tolerant scheme is developed for deterministic constraints. Adaptive fault-tolerant control scheme have been presented in \cite{Jin2015,Sun2020}, which can address deterministic constraints on state and input of the system. An adaptive fault-tolerant constrained control scheme has been developed to for commercial aircraft with actuator faults and constraints in \cite{Liu2019}. Fault-tolerant control of Euler-Lagrange systems has been discussed in \cite{Zhang2020}, where the output of the systems has to satisfy a deterministic constraint.} Recently, an RG-based reconfiguration scheme has been introduced in \cite{Li2021}. Even though the proposed scheme can effectively recover stability and constraint satisfaction properties after detecting the fault, it does not address the fault detection as it assumes that the fault can be detected immediately upon occurrence.

\begin{figure}
\centering
\includegraphics[width=7.5cm]{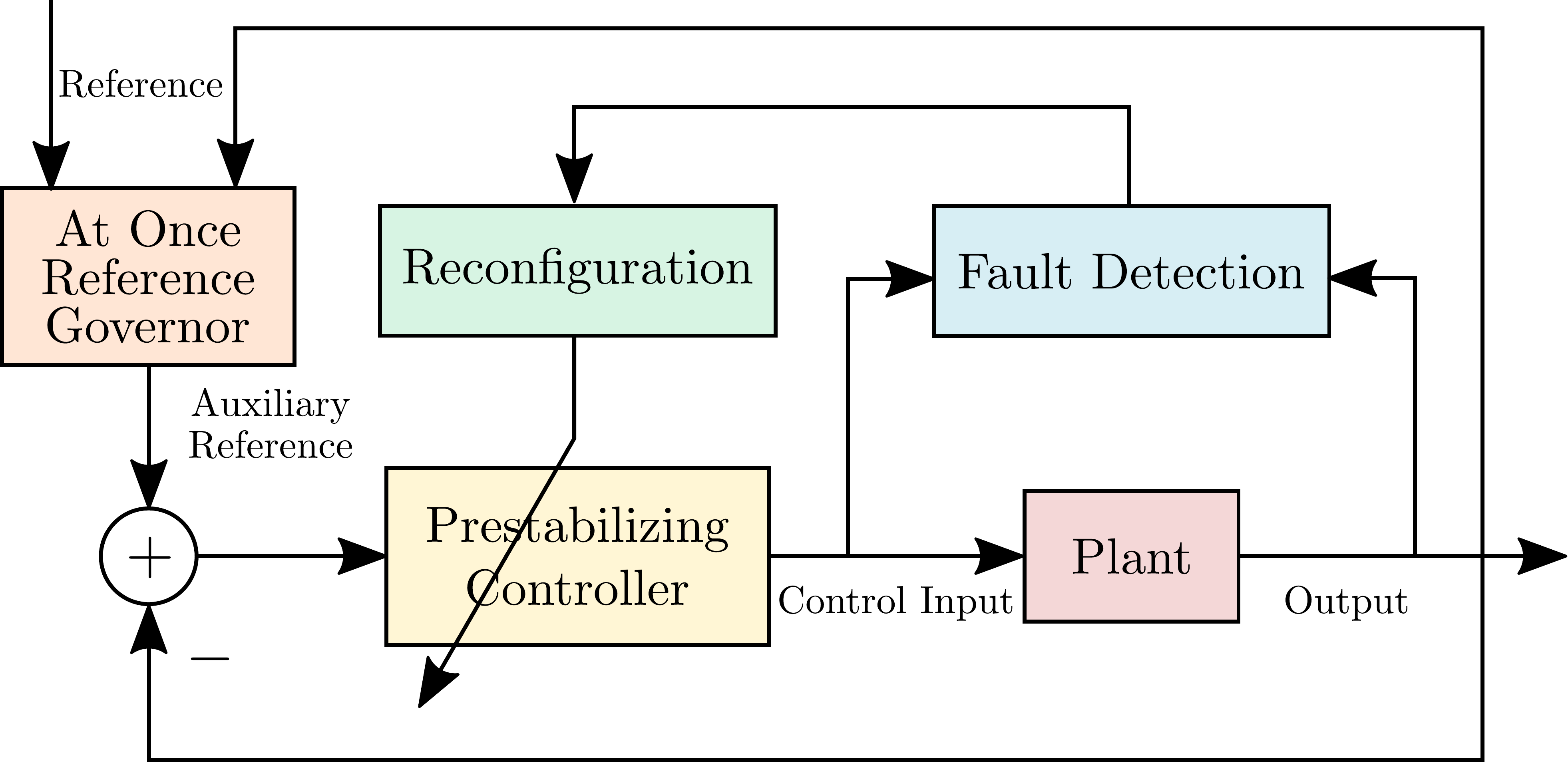}
\caption{Structure of the proposed Reference Governor-Based Fault-Tolerant  Constrained Control Scheme.}
\label{fig:structure}
\end{figure}

This paper proposes a RG-based fault-tolerant constrained control scheme, which addresses control, fault detection, and reconfiguration objectives. The structure of the proposed scheme is depicted in Fig. \ref{fig:structure}. {\color{black}Our motivation to use RG-based schemes is that they provide \textit{add-on} solutions, which can be attractive to practitioners interested in preserving an existing/legacy controller or concerned with computational burden and tuning complexity. Additionally, and as illustrated in this paper, they can non-conservatively restrict the operation of the system, which facilitates the ability of the system to recover from faults.} First, we propose At Once Reference Governor (AORG), which can be utilized to address tracking and constraint satisfaction requirements. This AORG is distinguished from the conventional RG by optimizing and applying the AR sequence over a time interval rather than at a given time instant. In order to detect the fault occurrence, we adopt the Multi-Model Adaptive Estimator (MMAE) \cite{Fekri2004,Hassani2009}. It will be shown that a bound on the detection performance can be expressed as an explicit function of the AR sequence over an interval. Two optimization problems will be formulated to determine the AR sequence during transient and at steady-state, such that to optimize the performance of the MMAE, while ensuring constraint satisfaction at all times. Finally, a reconfiguration scheme will be proposed to maintain functionality of the system despite the presence of the fault. This reconfiguration scheme is based on the recoverable sets \cite{McDonough2017,Li2021}.

The main contributions of this paper are: 1) presenting the AORG and proving its convergence and constraint-handling properties, 2) proving that the AORG can be integrated with the MMAE such that both control and detection objectives are addressed simultaneously, and 3) proposing a reconfiguration scheme to maintain the viability of the system despite the presence of the fault.

The remainder of this paper is organized as follows. Section \ref{sec:ProblemStatement} states the problem. Section \ref{sec:AORG} introduces the AORG scheme, and proves its constraint enforcement and convergence properties. In Section \ref{sec:FTC}, first, the MMAE is briefly discussed. It is then proven that, combined with the AORG, it is possible to determine the AR sequence such that the performance of the MMAE is optimized. A reconfiguration scheme is also proposed to recover system stability and constraint satisfaction property. Section \ref{sec:simulation} evaluates the effectiveness of the proposed scheme on a Boeing 747-100 aircraft in the presence of vertical stabilizer failure. Finally, Section \ref{sec:conclusion} concludes the paper.

%%%%%%%%%%%%%%%%%%%%%%%%%%%%%%%%%%%%%%%%%%%%%%%%%%%%%%%%%%%%%%%%%%%%%%%%%%%%%%%%
\paragraph*{Notation.}
$\mathbb{R}$ denotes the set of real numbers, and $\mathbb{R}_{\geq a}$ and $\mathbb{Z}_{\geq a}$ denote the real numbers and integer numbers greater than or equal to $a$, respectively. We denote the transpose of the matrix $R$ by $R^\top$. The Euclidean norm of a vector $x\in\mathbb{R}^n$ is denoted by $\left\Vert x\right\Vert=\sqrt{x_1^2+\cdots+x_n^2}$, whereas $\left\Vert x\right\Vert_R^2$ with $R=R^\top>0$ denotes the quadratic form $x^\top Rx$. The function $\mathcal{F}$ is used to represent the $\chi^2$ cumulative distribution function. For given sets $X,Y\subset\mathbb{R}^n$, $X\sim Y:=\{x:x+y\in X,\forall y\in Y\}$ is the Pontryagin set difference, and $X\oplus Y:\{x+y:x\in X,y\in Y\}$ is the Minkowski set sum. We use $|X|$ to represent the cardinality of the set $X$. The expected value of a random variable $x$ is denoted by $\mathbb{E}[x]$, and $\mathbb{P}(x\in U)$ indicates the probability that $x$ belongs to a certain event $U$. We denote the $n\times n$ identity matrix by $I_n$.  We denote the determinant function by $\text{det}(\cdot)$. $\mathcal{N}(\mu,\Sigma)$ indicates the Gaussian distribution with mean $\mu$ and covariance matrix $\Sigma$.

%%%%%%%%%%%%%%%%%%%%%%%%%%%%%%%%%%%%%%%%%%%%%%%%%%%%%%%%%%%%%%%%%%%%%%%%%%%%%%%%
\section{Problem Statement}\label{sec:ProblemStatement}
Consider the following discrete-time LTI system with multiple operating modes:
\begin{eqnarray}\label{eq:modesopenloop}
\left\{
\begin{array}{r@{=\,}l}
    x(t+1|\mu) & A_{o,\mu} x(t|\mu)+B_{o,\mu} u(t{\color{black}|\mu})+\omega_\mu(t) \\
    y(t|\mu) & Cx(t|\mu)+\xi(t)
\end{array}
\right.,
\end{eqnarray}
where $\mu\in\mathcal{M}=\{\mu_1,\cdots,\mu_{f}\}$ is the index of the mode of the system, $x(t|\mu)$ is the state of the system operating in mode $\mu$ at time $t$, $u(t)\in\mathbb{R}^p$ is the control input at time $t$, $y(t|\mu)$ is the output of the system operating in mode $\mu$ at time $t$, and the process noise $\omega_\mu(t)\in\mathbb{R}^n$ and the measurement noise $\xi(t)\in\mathbb{R}^m$ are mutually independent Gaussian processes with zero mean and covariance matrices $H_{\omega_\mu}\in\mathbb{R}^{n\times n}$ and $H_{\xi}\in\mathbb{R}^{m\times m}$, respectively. The $A_{o,\mu}$ and $B_{o,\mu}$ are the open loop dynamics and inputs matrices in operating mode $\mu$. The model \eqref{eq:modesopenloop} can represent a system with $f-1$ fault scenarios that manifest themselves as changes in the system matrices. The model $\mu_1$ corresponds to the nominal system operating without faults. {\color{black}Note that we assume that sensors are not affected by the faults.}

\begin{rem}
As an example, actuator faults can be captured by the model \eqref{eq:modesopenloop}.  In this paper, the failure of the $i$th actuator is represented by zeroing out the $i$th column of the matrix $B_o$.
\end{rem}

In order to stabilize the system in each mode, we use the following feedback plus feedforward control law:
\begin{eqnarray}\label{eq:controllaw}
u(t{\color{black}|\mu})=K_\mu x(t{\color{black}|\mu})+G_\mu v(t),
\end{eqnarray}
where $K_\mu\in\mathbb{R}^{p\times n}$ is the feedback gain matrix, $G_\mu\in\mathbb{R}^{p\times m}$ is the feedforward gain matrix, and $v(t)\in\mathbb{R}^m$ is the vector of reference commands (set-points). Thus, the closed-loop system takes the following form:
\begin{eqnarray}\label{eq:closedloopsystem}
\left\{
\begin{array}{r@{=\,}l}
    x(t+1|\mu) & A_\mu x(t|\mu)+B_\mu v(t)+\omega_\mu(t) \\
    y(t|\mu) & Cx(t|\mu)+\xi(t) 
\end{array}
\right.,
\end{eqnarray}
where $A_\mu=A_{o,\mu}+B_{o,\mu}K_\mu$ and $B_\mu=B_{o,\mu}G_\mu$. We assume that the feedback gain $K_\mu$ is such that $A_\mu$ is strictly Schur for all $\mu\in\mathcal{M}$.

Suppose that when the system is operating in mode $\mu$, its state and reference have to satisfy constraints of the following form,
\begin{eqnarray}\label{eq:constraints}
\left\{
\begin{array}{l}
\mathbb{E}\left[z_1(t|\mu)\right]\in\mathcal{Z}_1 \\
\mathbb{P}\left(z_2(t|\mu)\in\mathcal{Z}_2\right)\geq\beta
\end{array}
\right.,
\end{eqnarray}
where
\begin{eqnarray}
z_1(t|\mu)&=&L_xx(t|\mu)+L_vv(t)+\zeta(t),\\
z_2(t|\mu)&=&F_{x}x(t|\mu)+F_{v}v(t)+\varsigma(t),
\end{eqnarray}
are specified outputs, and where $L_x\in\mathbb{R}^{n_e\times n}$, $L_v\in\mathbb{R}^{n_e\times p}$, $n_e$ is the number of expectation constraints, $F_{x}\in\mathbb{R}^{n_c\times n}$, $F_{v}\in\mathbb{R}^{n_c\times p}$, $n_c$ is the number of chance constraints, $\zeta(t)\sim\mathcal{N}(0,H_{\zeta})$, and $\varsigma(t)\sim\mathcal{N}(0,H_{\varsigma})$. The constraint sets $\mathcal{Z}_1\subset\mathbb{R}^{n_e}$ and $\mathcal{Z}_2\subset\mathbb{R}^{n_c}$ are compact, convex, and contain the origin. {\color{black}Note that $z(t|\mu)=[z_1(t|\mu)~z_2(t|\mu)]^\top$ can represent real measurable outputs of the system, or can be used to represent the constraints on the state and input of the system.}

\begin{assum}\label{assumption:observable}
We assume that the pairs $(L_x,A_\mu)$ and $(F_{x},A_\mu)$ are observable for all $\mu$. 
\end{assum}

In this paper we consider the following problem:
\begin{prob}\label{prob:mainprob}
Consider system \eqref{eq:closedloopsystem} which is subject to constraints \eqref{eq:constraints}. Suppose that the operating mode of the system may change at anytime as a result of a fault/failure. Let $r(t)\in\mathbb{R}^m$ be the desired reference which is known over a sufficiently long preview window into the future at each time instant. For a given initial condition $x(0)$, develop a scheme to generate the AR signal $v(t)$ such that constraints \eqref{eq:constraints} are satisfied at all times, and for a constant desired reference, $v(t)$ tends to the desired reference $r(t)$. 
\end{prob}

To address this problem, we will first present the AORG scheme. Then, it will be shown that the AORG can be integrated with a fault detection scheme and a reconfiguration strategy, such that constraint satisfaction and convergence properties are retained.

%%%%%%%%%%%%%%%%%%%%%%%%%%%%%%%%%%%
\section{At Once Reference Governor}\label{sec:AORG}
{\color{black}As shown in Fig. \ref{fig:structure}, in the proposed fault-tolerant constrained control scheme, AORG is utilized to manipulate the AR sequence.} Unlike the conventional RG, AORG computes the AR sequence over a specified horizon. Unlike the conventional RG, AORG does not recompute the AR sequence at the next time instant rather it applies the AR sequence over this specified horizon and recomputes it at the end of this horizon for the next horizon.  Such an implementation is advantageous as it is able to improve fault detectability through reference manipulation {\color{black}(as will be shown in Subsection \ref{sec:ControlUnit})} and is applicable, e.g., to situations where commands represent waypoints which can be assumed to not change during each preview horizon. Note, however, that with AORG the system looses ability to respond to commands during each preview horizon.

{\color{black}In this section, we explain the general formulation of AORG, while its integration into the fault-tolerant constrained control scheme shown in Fig. \ref{fig:structure} will be discussed in Subsection \ref{sec:ControlUnit}. Consequently, and to simplify the notations, we drop the explicit dependence on $\mu$ in this section.} On the other hand, since we are concerned with the $k$ step ahead predictions made at time $t$, this will be reflected in the modified notations.

\subsection{Preliminaries}
{\color{black}In this subsection, we present two  propositions that will be used in this paper. The first proposition shows that the chance constraint given in \eqref{eq:constraints} can be enforced by enforcing a condition on the noise-free output. The second proposition shows how this can be done when the constraint set is a polytope.}

\begin{prop}[\cite{Du2018}]\label{prop:chanceconstraint}
Consider the following noise-driven and noise-free prediction models:
\begin{eqnarray}\label{eq:z1z2noisy}
\left\{
\begin{array}{r@{=\,}l}
\tilde{x}(k+1|t) & A\tilde{x}(k|t)+\omega(t+k),~\tilde{x}(0|t)=0 \\
\tilde{z}_1(k|t) & L_x\tilde{x}(k|t)+\zeta(t+k) \\
\tilde{z}_2(k|t) & F_{x}\tilde{x}(k|t)+\varsigma(t+k)
\end{array}
\right.,
\end{eqnarray}
and 
\begin{eqnarray}\label{eq:z1z2noisefree}
\left\{
\begin{array}{r@{=\,}l}
\hat{x}(k+1|t) & A\hat{x}(k|t)+Bv(t+k),~\hat{x}(0|t)=x(t)  \\
\hat{z}_1(k|t) & L_x\hat{x}(k|t)+L_vv(t+k) \\
\hat{z}_2(k|t) & F_x\hat{x}(k|t)+F_vv(t+k)
\end{array}
\right.,
\end{eqnarray}
where $k\in\mathbb{Z}_{\geq0}$, and note that $z_1(k|t)=\hat{z}_1(k|t)+\tilde{z}_1(k|t)$ and $z_2(k|t)=\hat{z}_2(k|t)+\tilde{z}_2(k|t)$. Then, $\mathbb{P}(z_2(t+k)\in\mathcal{Z}_2)\geq\beta$ if $\hat{z}_2(k|t)\oplus\mathcal{P}_\beta(k)\subset\mathcal{Z}_2$, where $\mathcal{P}_\beta(k)$ is the confidence ellipsoid with confidence level $\beta$ at time $k$, i.e., $\mathbb{P}(\tilde{z}_2(k|t)\in\mathcal{P}_\beta(k))=\beta$. %See Fig. \ref{fig:ChanceIllustration} for a geometric illustration. The $\mathcal{P}_\beta(k)$ can be computed as
%\begin{eqnarray}
%\mathcal{P}_\beta(k)=\left\{\theta\in\text{Im}(\Gamma(k)):\theta^\top\Gamma^+(k)\theta\leq\mathcal{F}^{-1}(\beta,m)\right\},
%\end{eqnarray}
%where $\Gamma(k)=F_{x}\Sigma(k)F_{x}^\top$, and the covariance matrix $\Sigma(k)$ can be obtained by the following covariance propagation equation:
%\begin{eqnarray}
%\Sigma(k+1)=A\Sigma(k)A^\top+H_{\omega},
%\end{eqnarray}
%with $\Sigma(0)$ as the covariance of the state measurement or estimate.
\end{prop}

\begin{prop}[\cite{Kolmanovsky1998}]
Suppose that
\begin{eqnarray}
\mathcal{Z}_2=\left\{z_2:z_{2,i}\leq \alpha_i,~i=1,\cdots,n_c\right\},
\end{eqnarray}
where $z_{2,i}$ is the $i$th element of $z_2$, and $\alpha_i\in\mathbb{R},~\forall i$. Then, \begin{eqnarray}\label{eq:constrainttightening}
\mathcal{Z}_2\sim\mathcal{P}_\beta(k)=\left\{z_2:z_{2,i}\leq \alpha_i-\sqrt{\mathcal{F}^{-1}(\beta,n_c)\Gamma_i(k)},~\forall i\right\},
\end{eqnarray}
where $\Gamma_i(k)$ is element $(i,i)$ of $\Gamma(k)=F_{x}\Sigma(k)F_{x}^\top+H_{\varsigma}$, with $\Sigma(k)$ obtained by the following covariance propagation equation,
\begin{eqnarray}\label{eq:SigmaFormulation}
\Sigma(k+1)=A\Sigma(k)A^\top+H_{\omega},
\end{eqnarray}
and where $\Sigma(0)$ is the covariance of the state measurement or estimate.
\end{prop}

%\begin{figure}
%\centering
%\includegraphics[width=6cm]{ChanceIllustration.eps}
%\caption{Caption}
%\label{fig:ChanceIllustration}
%\end{figure}

\subsection{AORG: Formulation and Properties}
Let $T>0$ be a chosen horizon and $[t,t+T]$ be a given command planning interval.

\subsubsection{The Maximal Output-Admissible Set}
The maximal output-admissible set is defined as the set of all initial states $x$ and input sequences $v_0,\cdots,v_T$, such that, assuming the input $v_T$ is constantly applied from the time instant $t+T$ onward, the ensuing outputs will always satisfy the constraints \eqref{eq:constraints}:
\begin{eqnarray}\label{eq:maximal2}
O_{\infty}&&=\Big\{(x,v_0,\cdots,v_T):\hat{z}_1(k|x,v_0,\cdots,v_T)\in\mathcal{Z}_1\text{ and }\nonumber\\
&&\hat{z}_2(k|x,v_0,\cdots,v_T)\in\mathcal{Z}_2\sim\mathcal{P}_\beta(k)\text{ for all } k\in\mathbb{Z}_{\geq0}\Big\},
\end{eqnarray}
where $\hat{z}_1(k|x,v_0,\cdots,v_T)$ and $\hat{z}_2(k|x,v_0,\cdots,v_T)$ for $k\in\mathbb{Z}_{\geq0}$ are\footnote{In the rest of this paper, $\hat{z}_1(\cdot)$ and $\hat{z}_2(\cdot)$ denote $\hat{z}_1(k|x,v_0,\cdots,v_T)$ and $\hat{z}_2(k|x,v_0,\cdots,v_T)$, respectively.}
\begin{eqnarray}
\hat{z}_1(\cdot)=&&L_xA^{k}\left(\hat{x}(T|x,v_0,\cdots,v_{T-1})-\left(I_n-A\right)^{-1}Bv_T\right)\nonumber\\
&&+\left(L_x(I_n-A)^{-1}B+L_v\right)v_T,\label{eq:z1hat1}\\
\hat{z}_2(\cdot)=&&F_xA^{k}\left(\hat{x}(T|x,v_0,\cdots,v_{T-1})-\left(I_n-A\right)^{-1}Bv_T\right)\nonumber\\
&&+\left(F_x(I_n-A)^{-1}B+F_v\right)v_T,\label{eq:z2hat1}
\end{eqnarray}
and $\hat{x}(T|x,v_0,\cdots,v_{T-1})$ can be computed via \eqref{eq:z1z2noisefree}, with the initial condition $\hat{x}(0|t)=x$, and $v(t+k)=v_k$ for $k\in\{0,\cdots,T-1\}$. Note that $O_\infty$ does not address constraints satisfaction within the interval $[t,t+T]$. Hence, in the following section, we will provide a method to compute a subset of $O_\infty$ which ensures constraints satisfaction at all times.

%$\hat{z}_1(k|x,v_0,\cdots,v_T)$ and $\hat{z}_2(k|x,v_0,\cdots,v_T),~\forall k\in\mathbb{Z}_{\geq0}$ can be computed via \eqref{eq:z1z2noisefree}, with initial condition $x$, and $v(k)=v_k$ for $k=0,\cdots,T-1$, and $v(k)=v_T$ for $k\in\mathbb{Z}_{\geq T}$. 

\subsubsection{A Constraint-Admissible Subset of $O_\infty$}
%Here, we provide a method to compute a close approximation to $O_\infty$. 
One possible way to compute a constrained-admissible subset of $O_\infty$ referred to 
as $\tilde{O}_\infty$, which ensures constraint satisfaction from time $t$ onward is to use the following set recursion:
\begin{eqnarray}\label{eq:recursiveOinf}
\tilde{O}_{k+1}=\tilde{O}_k\cap\Phi_{k+1},
\end{eqnarray}
where 
\begin{eqnarray}
\Phi_k=\big\{&&(x,v_0,\cdots,v_T):\hat{z}_1(k|x,v_0,\cdots,v_T)\in\mathcal{Z}_1\text{ and }\nonumber\\
&&\hat{z}_2(k|x,v_0,\cdots,v_T)\in\mathcal{Z}_2\sim\mathcal{P}_\beta(k)\big\},
\end{eqnarray}
with the initial condition $\tilde{O}_0=\Phi_0\cap\bar{\Phi}$, where
\begin{eqnarray}
\bar{\Phi}=\big\{&&(x,v_0,\cdots,v_T):\hat{z}_1(k|t)\in\mathcal{Z}_1\text{ and }\nonumber\\
&&\hat{z}_2(k|t)\in\mathcal{Z}_2\sim\mathcal{P}_\beta(k),~k=0,\cdots,T,\text{ and }\nonumber\\
&& \left(L_x(I_n-A)^{-1}B+L_v\right)v_T\oplus\mathcal{B}_\epsilon\subset\mathcal{Z}_1\text{ and }\nonumber\\
&&\left(F_x(I_n-A)^{-1}B+F_v\right)v_T\oplus\mathcal{B}_\epsilon\subset\bigcap\limits_{k=0}^{\infty}\left(\mathcal{Z}_2\sim\mathcal{P}_\beta(k)\right)\big\},\nonumber\\\label{eq:Phibar}
\end{eqnarray}
and $\hat{z}_1(k|t)$ and $\hat{z}_2(k|t)$ can be computed via \eqref{eq:z1z2noisefree}, with initial condition $\hat{x}(0|t)=x$, and $v(t+k)=v_k$ for $k=0,\cdots,T$, and $\mathcal{B}_\epsilon$ is an open ball of radius $\epsilon>0$. Simply, $\bar{\Phi}$ is the set of all initial conditions $x$ and input sequences $v_0,\cdots,v_T$ that steer the system such that the constraints are satisfied within the interval $[t,t+T]$, and $v_T$ is strictly steady-state admissible. Note that by Assumption \ref{assumption:observable}, $\bar{\Phi}$ is compact and convex. Thus, $\tilde{O}_0$ is compact and convex. {\color{black}Note that each recursive update in \eqref{eq:recursiveOinf} can be performed by simple offline linear algebra computations, and thus does not cause real-time implementation issues.}

%Consider the following set of steady-state admissible inputs:
%\begin{eqnarray}\label{eq:Oepsilon}
%O_{\epsilon}&&=\Big\{v_T:\left(L_x(I_n-A)^{-1}B+L_v\right)v_T\oplus\mathcal{B}_\epsilon\subset\mathcal{Z}_1,\nonumber\\
%&&\left(F_x(I_n-A)^{-1}B+F_v\right)v_T\oplus\mathcal{B}_\epsilon\subset\bigcap\limits_{k=0}^{\infty}\left(\mathcal{Z}_2\sim\mathcal{P}_\beta(k)\right)\Big\},
%\end{eqnarray}
%where $\mathcal{B}_\epsilon$ is an open ball of radius $\epsilon>0$. 

%The following lemma proves that $\tilde{O}_\infty$ computed via the set recursion given in \eqref{eq:recursiveOinf} is compact, convex, and finitely-determined. 

\begin{lem}\label{lemma:finitelydetermined}
The $\tilde{O}_\infty$ computed via the set recursion \eqref{eq:recursiveOinf} is compact and convex. Furthermore, it is finitely determined, i.e., there exists $k^\ast\in\mathbb{Z}_{\geq0}$ such that $\tilde{O}_{k^\ast}=\tilde{O}_\infty$.
\end{lem}

\begin{pf}
First, note that: 1) since  $\mathcal{Z}_2$ is compact and convex, it can be shown \cite{Kolmanovsky1998} that $\mathcal{Z}_2\sim\mathcal{P}_\beta(k)$ is compact and convex for all $k\in\mathbb{Z}_{\geq0}$, and 2) according to the set recursion given in \eqref{eq:recursiveOinf}, we have $\tilde{O}_{k+1}\subset\tilde{O}_k,~\forall k\in\mathbb{Z}_{\geq0}$.

We know that $\tilde{O}_0$ is compact  and convex. According to \eqref{eq:recursiveOinf}, $\tilde{O}_{1}$ is equal to $\tilde{O}_{0}\cap\Phi_{1}$. Note that $\Phi_{1}$ is closed and convex, as $\mathcal{Z}_1$ and $\mathcal{Z}_2\sim\mathcal{P}_\beta(1)$ are compact and convex. This means that $\tilde{O}_{1}$ is compact and convex. Therefore, by induction, it can be proven that $\tilde{O}_{k},~\forall k\in\mathbb{Z}_{\geq0}$ is compact and convex.

%Now, we begin to prove that $\tilde{O}_\infty$ is finitely-determined. 

According to \eqref{eq:z1hat1}-\eqref{eq:z2hat1}, since $A$ is assumed to be strictly Schur, for any $0<\epsilon'<\epsilon$ with $\epsilon$ as defined above, there exist $k'>0$ such that for all $k\in\mathbb{Z}_{\geq k'}$ we have
\begin{eqnarray}
L_xA^{k}\left(\hat{x}(T|x,v_0,\cdots,v_{T-1})-\left(I_n-A\right)^{-1}Bv_T\right)\in&&\mathcal{B}_{\epsilon'},\\
F_xA^{k}\left(\hat{x}(T|x,v_0,\cdots,v_{T-1})-\left(I_n-A\right)^{-1}Bv_T\right)\in&&\mathcal{B}_{\epsilon'},
\end{eqnarray}
where $\mathcal{B}_{\epsilon'}$ is an open ball of radius $\epsilon'$. Thus, according to \eqref{eq:Phibar}, for all $(x,v_0,\cdots,v_T)\in\tilde{O}_{k'}$ (and thus $(x,v_0,\cdots,v_T)$ bounded) and $k\in\mathbb{Z}_{\geq k'}$, we have 
\begin{eqnarray}
\hat{z}_1(\cdot)&\in&\mathcal{B}_{\epsilon'}\oplus\left(L_x(I_n-A)^{-1}B+L_v\right)v_T\subset\mathcal{Z}_1,\label{eq:z1steady}\\ 
\hat{z}_2(\cdot)&\in&\mathcal{B}_{\epsilon'}\oplus\left(F_x(I_n-A)^{-1}B+F_v\right)v_T\subset\bigcap\limits_{i=0}^{\infty}\left(\mathcal{Z}_2\sim\mathcal{P}_\beta(i)\right)\nonumber\\
&&~~~~~~~~~~~~~~~~~~~~~~~~~~~~~~~~~~~~~~~~\;\subset\mathcal{Z}_2\sim\mathcal{P}_\beta(k), \label{eq:z2steady}
\end{eqnarray}
which means that $(x,v_0,\cdots,v_T)\in\Phi_k$. Thus, according to \eqref{eq:recursiveOinf} and by induction, we have $(x,v_0,\cdots,v_T)\in\tilde{O}_{k}$. On the other hand, we know that $\tilde{O}_{k}\subset\tilde{O}_{k'},~\forall k\in\mathbb{Z}_{\geq k'}$. Therefore, there exists $k^\ast\leq k'$ such that  $\tilde{O}_{k^\ast}=\cdots=\tilde{O}_{k'}=\tilde{O}_{k'+1}=\cdots=\tilde{O}_\infty$. \qed
\end{pf}

%\begin{rem}
%Since the initial condition for the set recursion given in \eqref{eq:recursiveOinf} is compact and convex, despite conventional RGs, we do not set any assumption on the observability of the pairs $(L_x,A)$ and $(F_x,A)$.
%\end{rem}

%%%%%%%%%%%%%%%%%%%%%%%%%%%
\subsubsection{Determination of the AR}
Once $\tilde{O}_\infty$ is computed, the following AORG scheme can be employed to compute the AR within the interval $[t,t+T]$ by solving the following optimization problem: 
\begin{eqnarray}\label{eq:OPAOVRG}
\kappa_j^{i^\ast}=\left\{
\begin{array}{ll}
     &  \arg\;\max\limits_{\kappa_j^i,\forall i,j}\;\;\sum\limits_{i=0}^T\;\sum\limits_{j=0}^m\kappa_j^i\\
   \text{s.t.}  & \kappa_j^i\in[0,1],~j\in\{1,\cdots,m\}\\
   & v_0=v(t-1)+K_0(r(t)-v(t-1))\\
   & v_{i}=v_{i-1}+K_{i}(r(t+i)-v_{i-1}),~i=1,\cdots,T\\
   & (x(t),v_0,\cdots,v_T)\in \tilde{O}_{\infty}
\end{array}
\right.
\end{eqnarray}
where $K_i=\text{diag}\{\kappa_1^i\cdots,\kappa_m^i\}$, and then computing AR as $v(t+i)=v(t+i-1)+K_i^\ast(r(t+i)-v(t+i-1)),~i\in\{0,\cdots,T\}$. %Note that we can build an At Once Scalar Reference Governor (AOSRG) by setting $\kappa_1^i=\cdots=\kappa_m^i,~\forall i$.

\subsubsection{Infeasibility-Handling Mechanism}\label{rem:infeasibilityhandling}
Due to system \eqref{eq:closedloopsystem} having stochastic noise inputs, the computed AR by an AORG may not be recursively feasible. More specifically, the previously admissible reference $v(t+T)$ may be no longer constraint admissible at the beginning of the interval $[t+T+1,t+2T+1]$, i.e., $(x(t+T+1),v(t+T))\not\in\text{Proj}_{(x,v_0)}\tilde{O}_\infty$. In this case, as an infeasibility-handling mechanism, the reference will be kept unchanged for one step, i.e., $v(t+T+1)=v(t+T)$. The feasibility will be checked again at $t+T+1$. If it is feasible, the AORG will compute the AR over the shifted interval $[t+T+2,t+2T+2]$; or else, the AR will be kept unchanged for one more step.

\subsubsection{Properties}
The constraint-handling and convergence properties of the AORG will be proven in the following theorems. In order to prove some of these properties, we will follow a similar procedure to that of \cite{Kalabic2019}.

\begin{thm}\label{theorem:constrainhandling}
Consider the sequential distinct intervals with length of $T+1$, starting from $0$. Suppose that the AORG as in \eqref{eq:OPAOVRG} is used to compute the AR over the intervals. Also, suppose that the the infeasibility-handling mechanism described in Section \ref{rem:infeasibilityhandling} is employed. Then, constraints \eqref{eq:constraints} are satisfied at all times.
\end{thm}

\begin{pf}
Suppose that $(x(t),v(t-1))\in\text{Proj}_{(x,v_0)}\tilde{O}_\infty$ for some $t\in\mathbb{Z}_{\geq0}$. Let $k'=\inf\{k\in\mathbb{Z}_{\geq0}:v(t+K)\neq v(t-1)\}$, which is greater than 0 and possibly unbounded. Namely, $k'$ is the first instant after $t$ that the AR changes.

Regarding the expectation constraint, according to \eqref{eq:z1z2noisy} and \eqref{eq:z1z2noisefree}, and since $\tilde{O}_\infty\subset\Phi_k,~\forall k$, we have 
\begin{eqnarray}
\mathbb{E}\left[z_1(t+k)\right]=\mathbb{E}\left[\hat{z}_1(k|t)+\tilde{z}_1(k|t)\right]=\hat{z}_1(k|t)\in\mathcal{Z}_1,
\end{eqnarray}
for $k\in\{0,\cdots,k'-1\}$. Similarly, regarding the chance constraint, we have
\begin{eqnarray}
\mathbb{P}\left(z_2(t+k)\in\mathcal{Z}_2\right)&=&\mathbb{P}\left(\hat{z}_2(k|t)+\tilde{z}_2(k|t)\in\mathcal{Z}_2\right)\nonumber\\
&\geq&\mathbb{P}\left(\tilde{z}_2(k|t)\in\mathcal{P}_\beta(k)\right)=\beta,
\end{eqnarray}
for $k\in\{0,\cdots,k'-1\}$. Thus, the use of the infeasibility-handling mechanism described in Section \ref{rem:infeasibilityhandling} (i.e., the AR at the beginning of each interval changes only if the previous value is feasible) implies that $(x(t+k'),v(t+k'-1))\in\text{Proj}_{(x,v_0)}\tilde{O}_\infty$.

Therefore, assuming that $(x(0),v(-1))\in\text{Proj}_{(x,v_0)}\tilde{O}_\infty$ and consequently $(x(1),v(0))\in\text{Proj}_{(x,v_0)}\tilde{O}_\infty$, where $v(0)$ is computed by the AORG as in \eqref{eq:OPAOVRG}, by induction, it can be shown that the constraints \eqref{eq:constraints} are satisfied at all times. \qed
\end{pf}

\begin{thm}\label{theorem:eventualchange}
Suppose that at time $t$ which is the beginning of an interval, the previously admissible AR is no longer constraint admissible (i.e., $(x(t),v(t-1))\not\in\text{Proj}_{(x,v_0)}\tilde{O}_\infty$) and thus it remain unchanged. There exists a $k''\in\mathbb{Z}_{\geq0}$ such that $(x(t+k''),v(t-1))\in\text{Proj}_{(x,v_0)}\tilde{O}_\infty$, i.e., the system will eventually enter a configuration where it is safe to change the AR. 
\end{thm}

\begin{pf}
Suppose that $(x(t),v(t-1))\not\in\text{Proj}_{(x,v_0)}\tilde{O}_\infty$ for some $t\in\mathbb{Z}_{\geq0}$ which is the beginning time of an interval. This means that the infeasibility-handling mechanism discussed in Section \ref{rem:infeasibilityhandling} will keep the AR unchanged until the time that it is safe to change. Let $e(t+k)=x(t+k)-\hat{x}(k|t)$ be the prediction error at time $t+k,~\forall k\in\mathbb{Z}_{\geq0}$, where $\hat{x}(k|t)$ is as in \eqref{eq:z1z2noisefree}.

Since $A$ is strictly Schur, for any $0<\epsilon'<\epsilon/2$, there exists $k'$ such that for all $k\geq k'$ we have
\begin{eqnarray}
L_xA^{k}\left(x(t)-\left(I_n-A\right)^{-1}Bv(t-1)\right)\in&&\mathcal{B}_{\epsilon'},\\
F_xA^{k}\left(x(t)-\left(I_n-A\right)^{-1}Bv(t-1)\right)\in&&\mathcal{B}_{\epsilon'},
\end{eqnarray}
where $\mathcal{B}_{\epsilon'}$ is an open ball of radius $\epsilon'$. Furthermore, according to the Ergodic Theorem \cite{Coudene2016} and by defining a proper Gaussian measure \cite{Bogachev1998}, it can be shown \cite{Kalabic2019} that almost surely there exits $k''\geq k'$ such that $L_xA^{k}e(t+k''),F_xA^{k}e(t+k'')\in\mathcal{B}_{\epsilon'}$ for $k\in\{0,\cdots,k^\ast\}$, with $k^\ast$ as in Lemma \ref{lemma:finitelydetermined}. Thus, according to \eqref{eq:z1z2noisefree}, at time $t+k''$ and for $k\in\{0,\cdots,k^\ast\}$, we have
\begin{eqnarray}
\hat{z}_1(k|t+k'')=&&L_xA^{k}\left(x(t+k'')-\left(I_n-A\right)^{-1}Bv(t-1)\right)\nonumber\\
&&+\left(L_x(I_n-A)^{-1}B+L_v\right)v(t-1)\nonumber\\
=&&L_xA^{k+k''}\left(x(t)-\left(I_n-A\right)^{-1}Bv(t-1)\right)\nonumber\\
&&+\left(L_x(I_n-A)^{-1}B+L_v\right)v(t-1)\nonumber\\
&&+L_xA^ke(t+k''),\label{eq:z1kprime}
\end{eqnarray}
and similarly we have 
\begin{eqnarray}
\hat{z}_2(k|t+k'')=&&F_xA^{k+k''}\left(x(t)-\left(I_n-A\right)^{-1}Bv(t-1)\right)\nonumber\\
&&+\left(F_x(I_n-A)^{-1}B+F_v\right)v(t-1)\nonumber\\
&&+F_xA^ke(t+k'').\label{eq:z2kprime}
\end{eqnarray}

Equations \eqref{eq:z1kprime} and \eqref{eq:z2kprime}, together with \eqref{eq:z1steady} and \eqref{eq:z2steady}, imply that
\begin{eqnarray}
\hat{z}_1(k|t+k'')&\in&\mathcal{B}_{\epsilon'}\oplus\left(\mathcal{Z}_1\sim\mathcal{B}_\epsilon\right)\oplus\mathcal{B}_{\epsilon'}\nonumber\\
&\in&\mathcal{Z}_1\sim\mathcal{B}_{\epsilon-2\epsilon'}\subset\mathcal{Z}_1,\\
\hat{z}_2(k|t+k'')&\in&\mathcal{B}_{\epsilon'}\oplus\left(\bigcap\limits_{i=0}^{\infty}\left(\mathcal{Z}_2\sim\mathcal{P}_\beta(i)\right)\sim\mathcal{B}_\epsilon\right)\oplus\mathcal{B}_{\epsilon'}\nonumber\\
&\in&\bigcap\limits_{i=0}^{\infty}\left(\mathcal{Z}_2\sim\mathcal{P}_\beta(i)\right)\sim\mathcal{B}_{\epsilon-2\epsilon'}\subset\mathcal{Z}_2-\mathcal{P}_\beta(k),
\end{eqnarray}
which means that $(x(t+k''),v(t-1))\in\text{Proj}_{(x,v_0)}\tilde{O}_\infty$, and thus it is safe to change the AR. \qed
\end{pf}

\begin{thm}\label{theorem:asymptociconvergence}
Consider the sequential distinct intervals with length of $T+1$, starting from $0$. Suppose that $r(t)=r$, where $r$ is steady-state admissible\footnote{When $r$ is steady-state admissible, it means that $\left(L_x(I_n-A)^{-1}B+L_v\right)r\oplus\mathcal{B}_\epsilon\subset\mathcal{Z}_1$ and $\left(F_x(I_n-A)^{-1}B+F_v\right)r\oplus\mathcal{B}_\epsilon\subset\bigcap\limits_{k=0}^{\infty}\left(\mathcal{Z}_2\sim\mathcal{P}_\beta(k)\right)$, where $\mathcal{B}_\epsilon$ is an open ball of radius $\epsilon>0$.}. Then, $v(t)$ computed by the AORG as in \eqref{eq:OPAOVRG} asymptotically converges to $r$. 
\end{thm}

\begin{pf}
Suppose that $v(t)$ is the AR at time $t\in\mathbb{Z}_{\geq0}$. The AORG ensures that the convergence error, defined as the distance between the desired reference and the AR, is non-increasing. More precisely, the AORG ensures that $\left\Vert r-v(t+k)\right\Vert=\left\Vert r-v(t+k-1)\right\Vert$ if $v(t+k)=v(t+k-1)$, and  $\left\Vert r-v(t+k)\right\Vert <\left\Vert r-v(t+k-1)\right\Vert$ if $v(t+k)\neq v(t+k-1),~\forall k\in\mathbb{Z}_{\geq0}$. Thus, we only need to prove that if $v(t+k)=v(t+k-1)$ for a $k\in\mathbb{Z}_{\geq0}$, there exists a $k''>k$ such that $v(t+k'')\neq v(t+k-1)$.

Suppose that $v(t+k)=v(t+k-1)$. There are two reasons for this: 1) it is imposed by the infeasibility-handling scheme, and 2) it is the optimal solution obtained by \eqref{eq:OPAOVRG}. As proven in Theorem \ref{theorem:eventualchange} for the first reason, and following the same procedure of the proof of Theorem \ref{theorem:eventualchange} for the second reason, for any $0<\epsilon'<\epsilon/2$ there exists $k''$ such that for $k\in\{0,\cdots,k^\ast\}$ we have
\begin{eqnarray}
\hat{z}_1(k|t+k'')&\in&\mathcal{Z}_1\sim\mathcal{B}_{\epsilon-2\epsilon'},\\
\hat{z}_2(k|t+k'')&\in&\bigcap\limits_{i=0}^{\infty}\left(\mathcal{Z}_2\sim\mathcal{P}_\beta(i)\right)\sim\mathcal{B}_{\epsilon-2\epsilon'}.
\end{eqnarray}

Thus, if we set $v(t+k'')=v(t+k-1)+\Delta v$, where $\Delta v$ is an adjustment satisfying $\left(L_x(I_n-A)^{-1}B+L_v\right)\Delta v\in\mathcal{B}_{\epsilon''}$ with $\epsilon''<\epsilon-2\epsilon'$, we will have $(x(t+k''),v(t+k-1)+\Delta v)\in\text{Proj}_{(x,v_0)}\tilde{O}_\infty$. Since \eqref{eq:OPAOVRG} is convex, such adjustment will be always achieved upon existence. \qed
\end{pf}

\begin{rem}\label{remark:AORGProperties}
The AORG solves an optimization problem which is larger than that of the conventional RG, and hence can be more computationally demanding, in particular, for a large $T$. However, since the AORG optimizes AR sequence over an interval, in general, it results in improved (faster) tracking. Furthermore, optimizing over an interval enables the AORG, if augmented with detection schemes, to improve detection performance, which will be discussed in the following section.
\end{rem}

%%%%%%%%%%%%%%%%%%%%%%%%%%%%%%%%
\section{Fault Detection and Reconfiguration}\label{sec:FTC}
{\color{black} Thus far, we have presented the general formulation of an AORG. It has been proven that the AORG guarantees constraint satisfaction at all times, while ensuring reference tracking. Following the structure depicted in Fig. \ref{fig:structure},} in this section, we integrate the AORG with a detection scheme to identify the operating mode of the system. We also propose a reconfiguration scheme to maintain viability of the system after recognizing the mode change due to fault occurrence.

\begin{assum}\label{assumption:timebetweenfaults}
The time between subsequent faults/failures is large, implying that only one mode change needs to be considered at a time.
\end{assum}

%%%%%%%%%%%%%%%%%%

\subsection{Detection Unit}\label{sec:DetectionUnit}
\begin{figure}[!t]
\centering
\includegraphics[width=8cm]{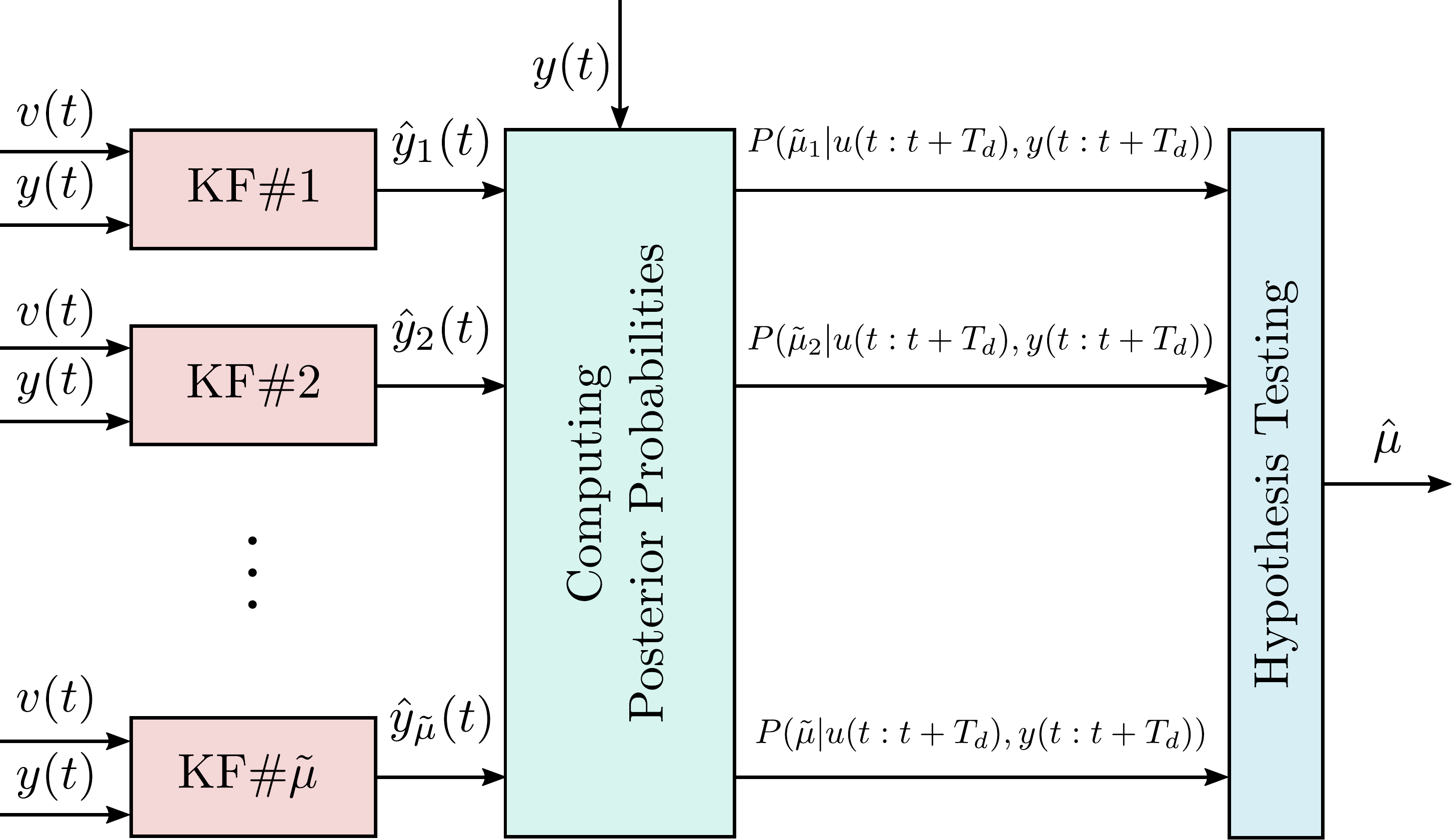}
\caption{The general structure of the MMAE deployed in this paper as the detection unit.}
\label{fig:KalmanBank}
\end{figure}

{\color{black}In this subsection, we show how AORG can be integrated with a detection unit to realize the structure shown in Fig. \ref{fig:structure}.} We employ the MMAE \cite{Hassani2009,Sadati2018} as the detection unit. Let $\mu$ be the current mode of the system, and $\mathcal{M}_{\mu}^+$ be the set of all successor modes of mode $\mu$. The MMAE involves the parallel operation of $|\mathcal{M}_\mu^+|+1$ Kalman filters, designed for systems $(A_{\bar{\mu}}+B_{\bar{\mu}}K_\mu,B_{\bar{\mu}}G_\mu,C,\textbf{0}),~\forall\bar{\mu}\in\{\mu\}\cup \mathcal{M}_\mu^+$, with $H_{\omega_{\bar{\mu}}}$ and $H_\xi$ as the process noise and measurement noise covariances. In MMAE, the residuals of the Kalman filters are used to identify the actual mode of the system. The general structure of the MMAE is shown in \figurename~\ref{fig:KalmanBank}, where $\hat{y}_i(t)$ is the predicted output of the $i$th Kalman filter.

\begin{rem}
In the MMAE, the actual mode can be identified (i.e., the posterior probabilities converge) almost surely \cite{Fekri2004,Rotondo2017}, if the systems are far apart. The distance between the systems can be assessed by means of Baram proximity metric \cite{Fekri2006} or gap metric \cite{Mahdianfar2011}. We assume that the Kalman filters are designed based upon system models that are sufficiently far apart. This assumption is reasonable, as the feedback and feedforward gains in \eqref{eq:controllaw} can be computed by optimizing the distance metrics.
\end{rem}

Let $\mathbb{P}\left(\tilde{\mu}|y(t:t+T_d),v(t:t+T_d-1)\right)$ be the posterior probability\footnote{Note that the Kalman filters are designed based upon systems that are not the ones defined in \eqref{eq:closedloopsystem}. To emphasize this difference, we use $\tilde{\mu}$ to denote the assumed modes in the design of Kalman filters. We also denote the set of these modes by $\tilde{\mathcal{M}}$.} of mode $\tilde{\mu}$ at time $t+T_d$ ($T_d\in\mathbb{Z}_{\geq0}$ is called detection time; see Remark \ref{rem:DetectionLength}) computed based on data over the time interval $[t,t+T_d]$, where $v(t:t+T_d-1):=[v(t)^\top,\cdots,v(t+T_d-1)^\top]^\top\in\mathbb{R}^{p(T_d+1)}$ and $y(t:t+T_d):=[y(t)^\top,\cdots,y(t+T_d)^\top]^\top\in\mathbb{R}^{m(T_d+1)}$. Note that $\mathbb{P}\left(\tilde{\mu}|y(t),v(t)\right)=\mathbb{P}(\tilde{\mu}),~\forall\tilde{\mu}$, where $\mathbb{P}(\tilde{\mu})$ is a known prior probability. Then the operating mode can be detected as the one which maximizes the posterior probability: 
\begin{eqnarray}\label{eq:detector}
\hat{\mu}=\arg\max\limits_{\tilde{\mu}}\mathbb{P}\left(\tilde{\mu}|y(t:t+T_d),v(t:t+T_d-1)\right).
\end{eqnarray}

%\begin{rem}\label{remark:consistency}
%Since the systems that the Kalman filters are design based upon are strictly proper, the input $v(t+T_d)$ has no effect on the posterior probability at time $t+T_d$. However, for the sake of consistency, we keep it in our formulations.
%\end{rem}

Let the detection objective be the quality of detection measured by the probability of mode misidentification by the detector \eqref{eq:detector}. This objective is a function of the AR sequence within the interval $[t,t+T_d]$, and can be expressed as
\begin{eqnarray}\label{eq:detectionobjectivefunction}
J_d(v(t:t+T_d-1))= \mathbb{E}\big[\sigma(\hat{\mu})\big],
\end{eqnarray}
where $\sigma(\hat{\mu})$ is zero when $\hat{\mu}$ determined by \eqref{eq:detector} is the actual operating mode of the system (i.e., the actual operating mode is identified correctly), and is 1 otherwise. The following theorem demonstrates that this objective function can be expressed as an explicit function of the AR sequence $v(t),\cdots,v(t+T_d-1)$.

\begin{thm}\label{theorem:detection}
Consider the MMAE shown in Fig. \ref{fig:KalmanBank}, and suppose that the AR sequence $v(t),\cdots,v(t+T_d-1)$ within the time interval $[t,t+T_d]$ has been specified. Then, the detection objective function \eqref{eq:detectionobjectivefunction} can be upper bounded by an explicitly computable function of the AR sequence.
\end{thm}

\begin{pf}
The detection objective function \eqref{eq:detectionobjectivefunction} can be expressed \cite{Puncochar2015,HosseinzadehACC2021} as
\begin{eqnarray}\label{eq:DetectionObjective1}
&&J_d(\cdot)=\mathbb{E}\big[\sigma(\hat{\mu})\big|v(t:t+T_d-1)\big]\nonumber\\
&&=\int\limits_{\mathbb{R}^{m(T_d+1)}}\sum\limits_{\tilde{\mu}\in\tilde{\mathcal{M}}}\sigma(\hat{\mu})\mathbb{P}\left(\tilde{\mu}\big|y(t:t+T_d),v(t:t+T_d-1)\right)\nonumber\\
&&~~~~~~~~~~\cdot \mathbb{P}\left(y(t:t+T_d)\big|v(t:t+T_d-1)\right)dy(t:t+T_d),
\end{eqnarray}

According to Bayes' theorem, since $0\leq\sigma(\hat{\mu})\leq1$, and due to the fact that the probability of the mode $\tilde{\mu}$ conditioned by only input data is equal to the \textit{a priori} probability of the mode $\tilde{\mu}$, \eqref{eq:DetectionObjective1} implies that
\begin{eqnarray}
J_d(\cdot)\leq\int\limits_{\mathbb{R}^{m(T_d+1)}}\sum\limits_{\tilde{\mu}\in\tilde{\mathcal{M}}}&&\mathbb{P}\left(y(t:t+T_d)\big|\tilde{\mu},v(t:t+T_d-1)\right)\nonumber\\
&&\cdot \mathbb{P}(\tilde{\mu})dy(t:t+T_d).\label{eq:Jd1}
\end{eqnarray}

Following the same arguments as in \cite{Blackmore2006}, the right-hand side of \eqref{eq:Jd1} can be upper bounded by $\hat{J}_d(v(t:t+T_d-1))$, which can be computed as
\begin{eqnarray}
\hat{J}_d(v(t:t+T_d))=\frac{1}{2}\sum\limits_{\tilde{\mu}\in\tilde{\mathcal{M}}}\sum\limits_{\check{\mu}\in\tilde{\mathcal{M}}}\sqrt{\mathbb{P}(\tilde{\mu})\mathbb{P}(\check{\mu})}e^{-\rho_{\tilde{\mu}\check{\mu}}},\label{eq:Jdhat}
\end{eqnarray}
where 
\begin{eqnarray}
\rho_{\tilde{\mu}\check{\mu}}=&&\frac{1}{4}\left(\eta_{\tilde{\mu}}-\eta_{\check{\mu}}\right)^\top\left(\Psi_{\tilde{\mu}}+\Psi_{\check{\mu}}\right)^{-1}\left(\eta_{\tilde{\mu}}-\eta_{\check{\mu}}\right)\nonumber\\
&&+\frac{1}{2}\ln\left(\frac{\text{det}\left(\frac{\Psi_{\tilde{\mu}}+\Psi_{\check{\mu}}}{2}\right)}{\sqrt{\text{det}(\Psi_{\tilde{\mu}})\text{det}(\Psi_{\check{\mu}})}}\right),\label{eq:phi}
\end{eqnarray}
with
\begin{eqnarray}
\eta_{\tilde{\mu}}=&&\left[\left(C\hat{x}(0|t,\tilde{\mu})\right)^\top~\cdots~\left(C\hat{x}(T_d|t,\tilde{\mu})\right)^\top\right]^\top,\\
\eta_{\check{\mu}}=&&\left[\left(C\hat{x}(0|t,{\check{\mu}})\right)^\top~\cdots~\left(C\hat{x}(T_d|t,{\check{\mu}})\right)^\top\right]^\top,
\end{eqnarray}
\begin{eqnarray}
\Psi_{\tilde{\mu}}=&&\Big[C\Sigma(0)C^\top~C\Sigma(1|\tilde{\mu})C^\top+H_\xi~~C\Sigma(2|\tilde{\mu})C^\top+H_\xi\nonumber\\
&&~~\cdots~C\Sigma(T_d|\tilde{\mu})C^\top+H_\xi\Big],\\
\Psi_{\check{\mu}}=&&\Big[C\Sigma(0)C^\top~C\Sigma(1|{\check{\mu}})C^\top+H_\xi~C\Sigma(1|{\check{\mu}})C^\top+H_\xi\nonumber\\
&&~~\cdots~C\Sigma(T_d|{\check{\mu}})C^\top+H_\xi\Big],
\end{eqnarray}
in which $\hat{x}(k|t,\tilde{\mu})$ and $\hat{x}(k|t,\check{\mu}),~k\in\{0,\cdots,T_d\}$, are as in \eqref{eq:z1z2noisefree}, and $\Sigma(k|\tilde{\mu})$ and $\Sigma(k|\check{\mu}),~k\in\{0,\cdots,T_d\}$, are as in \eqref{eq:SigmaFormulation}, computed with matrices of modes $\tilde{\mu}$ and $\check{\mu}$, respectively.

Note that $\eta_{\tilde{\mu}}$ and $\eta_{\check{\mu}}$ are explicit functions of the AR sequence $v(t),\cdots,v(t+T_d-1)$, and it can be easily shown \cite{HosseinzadehACC2021} that $\rho_{\tilde{\mu}\check{\mu}}$ is a quadratic function of this AR sequence. This completes the proof. \qed
% Thus, the objective function \eqref{eq:detectionobjectivefunction} can be upper bounded with an explicit quadratic function of the AR sequence $v(t),\cdots,v(t+T_d-1)$.\qed %Note that as discussed in Remark \ref{remark:consistency}, $v(t+T_d)$ has no effect on the objective function \eqref{eq:detectionobjectivefunction}.
\end{pf}

\begin{rem}\label{remark:convexityJ_d}
The upper-bound $\hat{J}_d(v(t:t+T_d-1))$ given in \eqref{eq:Jdhat} is, in general, non-convex in the AR sequence $v(t),\cdots,v(t+T_d-1)$.
\end{rem}

%%%%%%%%%%%%%%%%%%%%%%%%%%%%%%%%

\subsection{Reconfiguration Scheme}
{\color{black}As shown in Fig. \ref{fig:structure}, a reconfiguration scheme is employed in the proposed fault-tolerant constrained control scheme to bring the system to a safe configuration associated with the new operating mode, once a mode change is detected. In this subsection, we design the reconfiguration scheme.} Before proceeding, we make the following assumption.

\begin{assum}\label{assumption:T_r}
We assume that some of the constraints can be temporarily relaxed. This assumption is reasonable, as in practice, constraints are often imposed conservatively to extend system operating life \cite{Kolmanovsky2011,Li2021}. In mathematical terms, by relaxation we mean that $\mathbb{E}\left[z_1(t|\mu)\right]\in\mathcal{Z}_1^+$ and $\mathbb{P}\left(z_2(t|\mu)\in\mathcal{Z}_2^+\right)\geq\beta$, where  $\mathcal{Z}_1^+\supset\mathcal{Z}_1$ and $\mathcal{Z}_2^+\supset\mathcal{Z}_2$ are the extended sets. Also, we assume that this relaxation can last up to maximum $T_e\in\mathbb{Z}_{\geq2}$ time steps, which is referred to as the maximum extension time.
\end{assum}

Let $T_r$ ($0<T_r<T_e$) be the specified recovery time, i.e., the time within the system should complete the recovery and enter a safe configuration. Let $\mathcal{R}_{\mu}^{T_r}$ be a recoverable set associated with mode $\mu$. This set contains all states that can be steered into $\text{Proj}_{x}\tilde{O}_{\infty,\mu}$ within $T_r$ time steps, i.e.,
\begin{eqnarray}
\mathcal{R}_{\mu}^{T_r}=&&\big\{x(t): \exists~v(t),\cdots,v(t+T_r-1) \text{ such that }\nonumber\\
&&\hat{z}_1(k|t,\mu)\in\mathcal{Z}_1^+\text{ and }\hat{z}_2(k|t,\mu)\in\mathcal{Z}_2^+\sim\mathcal{P}_\beta(k),\nonumber\\
&&\text{for all }k\in\{0,\cdots,T_r-1\},\text{ and}\nonumber\\
&&(\hat{x}(T_r|t,\mu),v(t+T_r-1))\in\text{Proj}_{(x,v_0)}\tilde{O}_{\infty,\mu}\big\},\label{eq:recoverableset}
\end{eqnarray}
where $\hat{x}(k|t,\mu)$, $\hat{z}_1(k|t,\mu)$, and $\hat{z}_2(k|t,\mu)$ are as in \eqref{eq:z1z2noisefree}. {\color{black}Note that this recoverable set can be computed once and offline for all $\mu$.} 
%Note that there is no need to consider the effects of the random disturbances in the last line, as they are already considered in $\tilde{O}_{\infty,\mu}$. \textcolor{red}{See Fig. \ref{fig:RecoverableSet} for a geometric illustration.} %The following theorem proposes a reconfiguration scheme.

\begin{lem}\label{lemma:reconfiguration}
Suppose that at time $t$, a mode change from $\mu$ to $\bar{\mu}$ is detected, and $x(t)\in\mathcal{R}_{\bar{\mu}}^{T_r}$. Then, the AR sequence $v(t),\cdots,v(t+T_r-1)$ computed via the following optimization problem:
\begin{eqnarray}
\left\{
\begin{array}{ll}
&  \min\limits_{v(t),\cdots,v(t+T_r-1)}\;\sum\limits_{i=0}^{T_r-1}\left\Vert v(t+i)-r\right\Vert_R^2\\
\text{s.t.}  & \hat{z}_1(i|t,\bar{\mu})\in\mathcal{Z}_1^+\text{ and }\hat{z}_2(i|t,\bar{\mu})\in\mathcal{Z}_2^+\sim\mathcal{P}_\beta(i)\\
&\text{for all }i\in\{0,\cdots,T_r-1\}\\
& (\hat{x}(T_r|t,\bar{\mu}),v(t+T_r-1))\in\text{Proj}_{(x,v_0)}\tilde{O}_{\infty,\bar{\mu}}
\end{array}
\right.,
\end{eqnarray}
provides a safe reconfiguration. 
\end{lem}

\begin{pf}
The existence of such AR sequence $v(t),\cdots,v(t+T_r-1)$ follows from the definition of the set $\mathcal{R}_{\bar{\mu}}^{T_r}$ and Assumption \ref{assumption:timebetweenfaults}. The AORG scheme can then be utilized to control the system from time $t+T_r$ onward. \qed
\end{pf}

%%%%%%%%%%%%%%%%
\subsection{Control Unit}\label{sec:ControlUnit}
In this section, we discuss how to employ the AORG to generate the AR sequence $v(t),\cdots,v(t+T_d-1)$ such that the detection scheme presented in Section \ref{sec:DetectionUnit} detects the fault with high probability of correctness, without breaking the control objectives, while ensuring a safe recovery upon detection. %Before starting, we set the following remarks. See Fig. \ref{fig:AdmissibleSets} for a geometric illustration of the notions. 

\subsubsection{Control During Transient}\label{sec:controlduringtransient}
The following theorem addresses an active fault detection and control problem based upon the AORG scheme. %In the rest of this paper, $\tilde{O}_{\infty,\mu}$ and $\tilde{O}_{\infty,\mu}^+$ indicate the approximation to the maximal constraint admissible set of the system operating in mode $\mu$, created based upon the constraint sets $\mathcal{Z}_1$ and $\mathcal{Z}_2$, and $\mathcal{Z}_1^+$ and $\mathcal{Z}_2^+$, respectively.

\begin{thm}\label{theorem:activetransient}
Suppose that $\mu$ is the operating mode of the system at time $t$ which is the beginning of a detection interval. Let $(x(t),v(t-1))\in\text{Proj}_{(x,v_0)}\tilde{O}_{\infty,\mu}$, and $\left\Vert v(t-1)-r\right\Vert>\vartheta$ for some $\vartheta\in\mathbb{R}_{\geq0}$. Suppose that the AR sequence $v(t),\cdots,v(t+T_d-1)$ is computed via the optimization problem \eqref{eq:OPAOVRG} with the following cost function,
\begin{eqnarray}\label{eq:CostFunctionCD}
\Omega\sum\limits_{i=0}^{T_d-1}\sum\limits_{j=0}^{m}\kappa_j^i-(1-\Omega)\hat{J}_d(v_0,\cdots,v_{T_d-1}),
\end{eqnarray}
and with the following extra constraints:
\begin{eqnarray}
\hat{x}(T_d|t,\bar{\mu},K_\mu,G_\mu)&\in&\mathcal{R}_{\bar{\mu}}^{T_r},\label{eq:extraconstraint1}\\
\hat{z}_1(k|t,\bar{\mu},K_\mu,G_\mu)&\in&\mathcal{Z}_1^+,\label{eq:extraconstraint2}\\
\hat{z}_2(i|t,\bar{\mu},K_\mu,G_\mu)&\in&\mathcal{Z}_2^+\sim\mathcal{P}_\beta(k),\label{eq:extraconstraint3}
\end{eqnarray}
for all $k\in\{0,1,\cdots,T_d\}$ and all $\bar{\mu}\in\mathcal{M}_{\mu}^+$, where $\Omega\in[0,1]$ is a design parameter that defines the trade-off between control performance and detection quality, $\hat{x}(\cdot)$, $\hat{z}_1(\cdot)$, and $\hat{z}_2(\cdot)$ are as in \eqref{eq:z1z2noisefree},  $\hat{J}_d(\cdot)$ is as in \eqref{eq:Jdhat}, $\mathcal{R}_{\mu}^{T_r}$ is as in \eqref{eq:recoverableset}, $\mathcal{Z}_1^+$ and $\mathcal{Z}_2^+$ are the extended sets as discussed in Assumption \ref{assumption:T_r}, and $T_d$ is the detection time. Then, the AR sequence $v(t+i)=v_i^\ast,~i\in\{0,\cdots,T_d-1\}$, where $v_i^\ast$ is the optimal solution of the above problem, results in the relaxed constraints being satisfied within the interval $[t,t+T_d-1]$. If the operating mode of the system remains constant over the interval $[t,t+T_d-1]$, the computed AR sequence results in the constraints \eqref{eq:constraints} being satisfied within the above-mentioned interval. 
\end{thm}

\begin{pf}
The proof is a straightforward application of Theorems \ref{theorem:constrainhandling} and \ref{theorem:detection}, Assumptions \ref{assumption:timebetweenfaults} and \ref{assumption:T_r}, and Remarks \ref{rem:DetectionLength}.
\end{pf}

\begin{rem}
Another possible way to pursue both control and detection aims is to let one of the objective functions to take arbitrary value up to a known upper limit value, and then to enforce this as a constraint and minimize the other objective function. In this paper we only study the convex combination of two objective functions, as in \eqref{eq:CostFunctionCD}. 
\end{rem}

\begin{rem}\label{rem:DetectionLength}
As discussed in \cite{Fekri2004,Sadati2018}, the MMAE may not identify the mode change in one interval if it occurs at a time which is close to the end of the detection interval. According to this fact and Assumption \ref{assumption:T_r}, the following condition should hold true:
\begin{eqnarray}
T_r<T_e-2T_d,
\end{eqnarray}
where $T_r$ is the recovery time as in \eqref{eq:recoverableset}, $T_e$ is the extension time as in Assumption \ref{assumption:T_r}, and $T_d$ is the detection time. 
\end{rem}

\begin{figure}[!t]
\centering
\includegraphics[width=7.5cm]{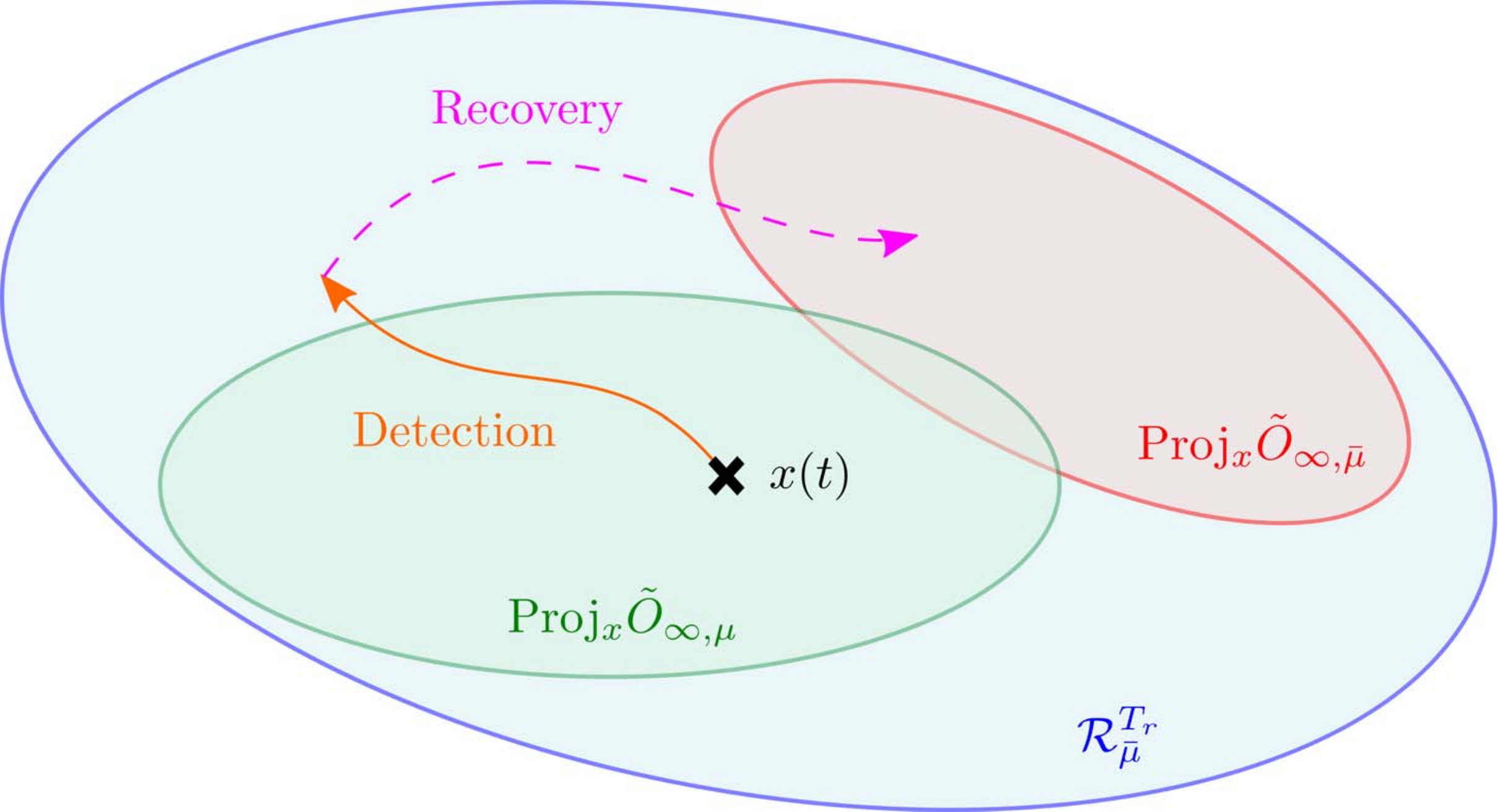}
\caption{Geometric illustration of the recoverable set and the maximal output-admissible set.}
\label{fig:RecoverableSet}
\end{figure}

\begin{rem}\label{remark:recoverablesetcondition}
According to the last constraint in \eqref{eq:OPAOVRG} and the constraint \eqref{eq:extraconstraint1}, to satisfy the tracking properties the following condition should hold true for all $\mu$ and for all $\bar{\mu}\in \mathcal{M}_{\mu}^+$:
\begin{eqnarray}\label{eq:conditionrecovery}
\text{Proj}_{x}\tilde{O}_{\infty,\mu}\subseteq\mathcal{R}_{\bar{\mu}}^{T_r}.
\end{eqnarray}
See \figurename~\ref{fig:RecoverableSet} for a geometric illustration. As discussed in \cite{Li2021}, there are three mutually non exclusive approaches to satisfy this condition: 1) to increase $T_r$, i.e., increase $T_e$ and/or decrease $T_d$; 2) to reduce $\mathcal{Z}_1$ and $\mathcal{Z}_2$ (i.e., to tighten constraints); and 3) to enlarge $\mathcal{Z}_1^+$ and $\mathcal{Z}_2^+$ (i.e., to further relax the constraints). Note that when the condition \eqref{eq:conditionrecovery} holds, if the state $x(t)$ does not belong to the recoverable set at the time of detection due to the random disturbances, the infeasibility-handling mechanism ensures that the state will eventually enter the recoverable set.
\end{rem}

\begin{figure}[!t]
\centering
\includegraphics[width=7.5cm]{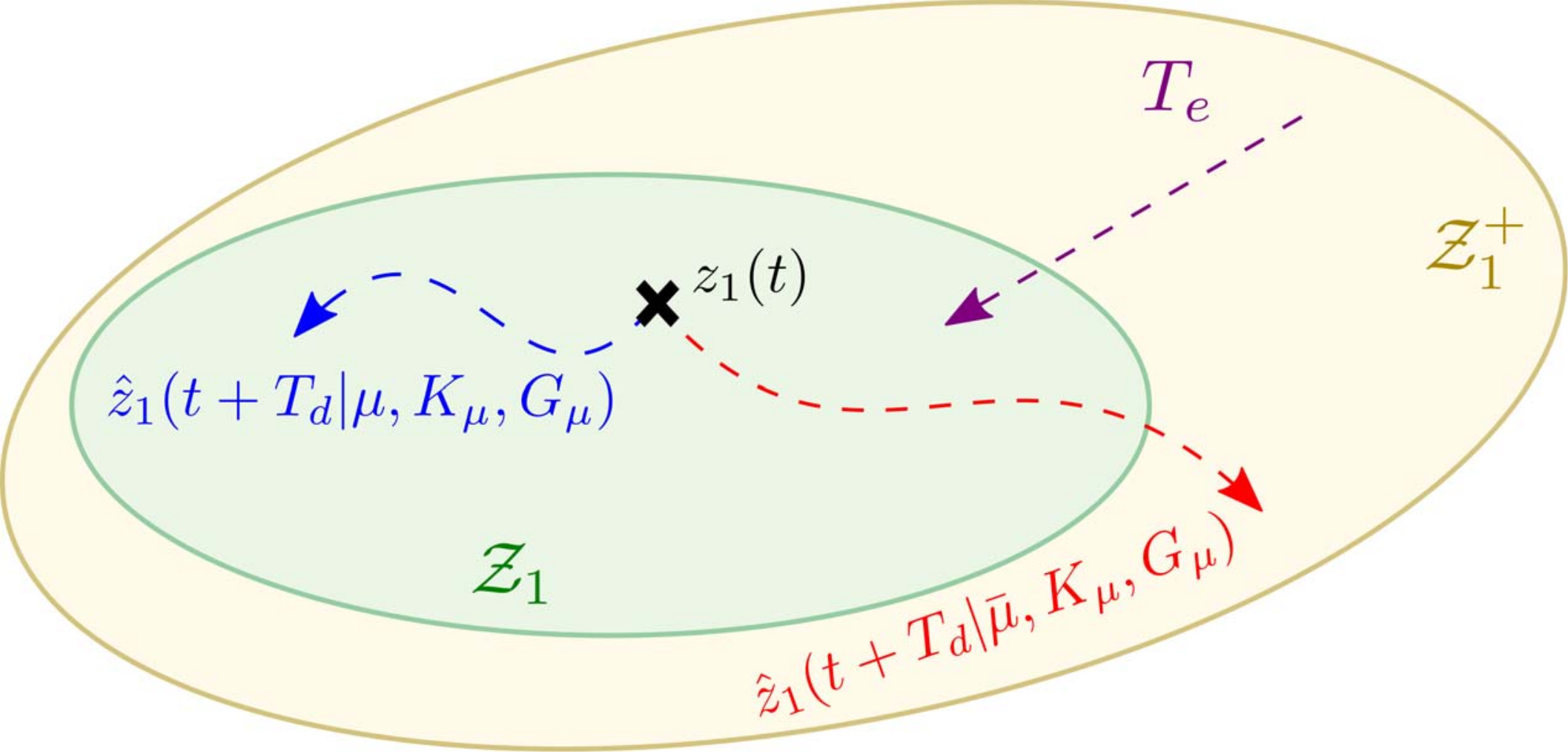}
\caption{Geometric illustration of the constraint admissible set $\mathcal{Z}_1$, extended admissible set $\mathcal{Z}_1^+$, maximum extension time $T_e$, and detection time $T_d$.}
\label{fig:AdmissibleSets}
\end{figure}

\begin{rem}\label{remark:T_v}
Constraints \eqref{eq:extraconstraint2} and \eqref{eq:extraconstraint3} mean that if the control law and the AR sequence designed for system $\mu$ are applied to system $\bar{\mu}$, $\hat{z}_1(\cdot)$ and $\hat{z}_2(\cdot)$ do not exit the extended sets $\mathcal{Z}_1^+$ and $\mathcal{Z}_2^+$, respectively. See Fig. \ref{fig:AdmissibleSets} for a geometric illustration (this figure shows only the set $\mathcal{Z}_1$). There are three mutually non exclusive approaches to satisfy these constraints (and consequently ensure recursive feasibility): 1) to reduce $T_d$; 2) to reduce $\mathcal{Z}_1$ and $\mathcal{Z}_2$; and 3) to enlarge $\mathcal{Z}_1^+$ and $\mathcal{Z}_2^+$.
\end{rem}

\begin{rem}\label{remark:guaranteemisidentification}
Due to the asymptotic convergence of the MMAE, correct identification may not be guaranteed for small $T_d$. Thus, the reconfiguration scheme may be applied based upon a misidentification. One intuitive way to cope with this issue is to apply the reconfiguration scheme if the mode change is being detected in two consecutive intervals.
\end{rem}

\begin{rem}
The optimization problem outlined in Theorem \ref{theorem:activetransient} is, in general, non-convex (see Remark \ref{remark:convexityJ_d}). However, its solution can be computed numerically by means of available nonlinear programming tools, e.g., \texttt{bmibnb} \cite{Lofberg2004} and \texttt{GloptiPoly3} \cite{Henrion2004}.
\end{rem}

\begin{rem}
Suppose that the operating mode of the system remains unchanged, and $r$ is steady-state admissible. Suppose that we employ the optimization problems mentioned in Theorem \ref{theorem:activetransient} in the following intervals, and we use the infeasibility-handling mechanism presented in Section \ref{rem:infeasibilityhandling}. Then, according to Theorems \ref{theorem:eventualchange}, \ref{theorem:asymptociconvergence}, and \ref{theorem:activetransient}, it can be shown that constraints \eqref{eq:constraints} are satisfied at all times, and $v(t)$ asymptotically converges to $r$. 
\end{rem}

\begin{rem}\label{remark:convergencedetection}
Suppose that a fault occurs at a time close to the end of the interval $[t,t+T_d]$. As mentioned in Remark \ref{rem:DetectionLength}, the MMAE may not be able to detect the fault at time $t+T_d$. If $(x(t+T_d),v(t+T_d-1))\in\text{Proj}_{(x,v_0)}\tilde{O}_{\infty,\mu}$, the fault will be detected at the end of the next interval, i.e., at time $t+2T_d$. Or else, since the infeasibility-handling mechanism presented in Section \ref{rem:infeasibilityhandling} keeps the AR, the posterior probabilities will go on evolving, and thus the fault will be detected in a few time steps ($<T_d$) \cite{HosseinzadehACC2021}. 
\end{rem}

\subsubsection{Control At Steady-State}
The following theorem formulates an active fault detection and control problem at steady-state. Note that discussions in Remark \ref{rem:DetectionLength}-\ref{remark:convergencedetection} hold true in this case.

\begin{thm}\label{theorem:activesteadystate}
Suppose that $\mu$ is the operating mode of the system at time $t$ which is the beginning of an interval. Let $(x(t),v(t-1))\in\text{Proj}_{(x,v_0)}\tilde{O}_{\infty,\mu}$, and $\left\Vert v(t-1)-r\right\Vert\leq\vartheta$. 
Consider the following optimization problem:
\begin{eqnarray}
\left\{
\begin{array}{ll}
     &  \min\limits_{v_0,\cdots,v_{T_d-1}}\;\Omega\sum\limits_{i=0}^{T_d-1}\left\Vert v_i-r\right\Vert_R^2\\
     &~~~~~~~~~~~~~~~~+(1-\Omega)\hat{J}_d(v_0,\cdots,v_{T_d-1})\\
   \text{s.t.} 
       & \left\Vert v_i-r\right\Vert\leq\vartheta,~\forall i\in\{0,\cdots,T_d-1\}\\
     & (x(t),v_0,\cdots,v_{T_d-1})\in \tilde{O}_{\infty,\mu}\\
   & \hat{x}(T_d|t,\bar{\mu},K_\mu,G_\mu)\in\mathcal{R}_{\bar{\mu}}^{T_r}\\
 & \hat{z}_1(k|t,\bar{\mu},K_\mu,G_\mu)\in\mathcal{Z}_1^+\\
& \hat{z}_2(k|t,\bar{\mu},K_\mu,G_\mu)\in\mathcal{Z}_2^+\sim\mathcal{P}_\beta(k)
\end{array}
\right.,
\end{eqnarray}
for $k\in\{0,1,\cdots,T_d\}$ and $\forall\bar{\mu}\in\mathcal{M}_{\mu}^+$, where $R=R^\top>0$ is a design matrix, $\Omega\in[0,1]$ is a design parameter. Then, the AR sequence $v(t+i)=v_i^\ast,~i\in\{0,\cdots,T_d-1\}$, where $v_i^\ast$ is the optimal solution, results in the relaxed constraints being satisfied within the interval $[t,t+T_d-1]$. If the operating mode of the system remains constant over the interval $[t,t+T_d-1]$, the computed AR sequence results in the constraints \eqref{eq:constraints} being satisfied within the above-mentioned interval.
\end{thm}

\begin{pf}
The proof is a straightforward application of Theorems \ref{theorem:constrainhandling} and \ref{theorem:detection}, Assumptions \ref{assumption:timebetweenfaults} and \ref{assumption:T_r}, and Remarks \ref{rem:DetectionLength}.
\end{pf}

\begin{rem}
The idea of Theorem \ref{theorem:activesteadystate} can be interpreted as injecting a small perturbation signal to the system in steady-state to improve detection performance. Note that injecting a perturbation signal for diagnostics purposes is exploited in many real-world systems (e.g., \cite{Tousi2011,Du2018})
\end{rem}

%%%%%%%%%%%%%%%%%%%
\section{Simulation Study}\label{sec:simulation}
\begin{figure}
    \centering
    \includegraphics[width=8cm,height=4.8cm]{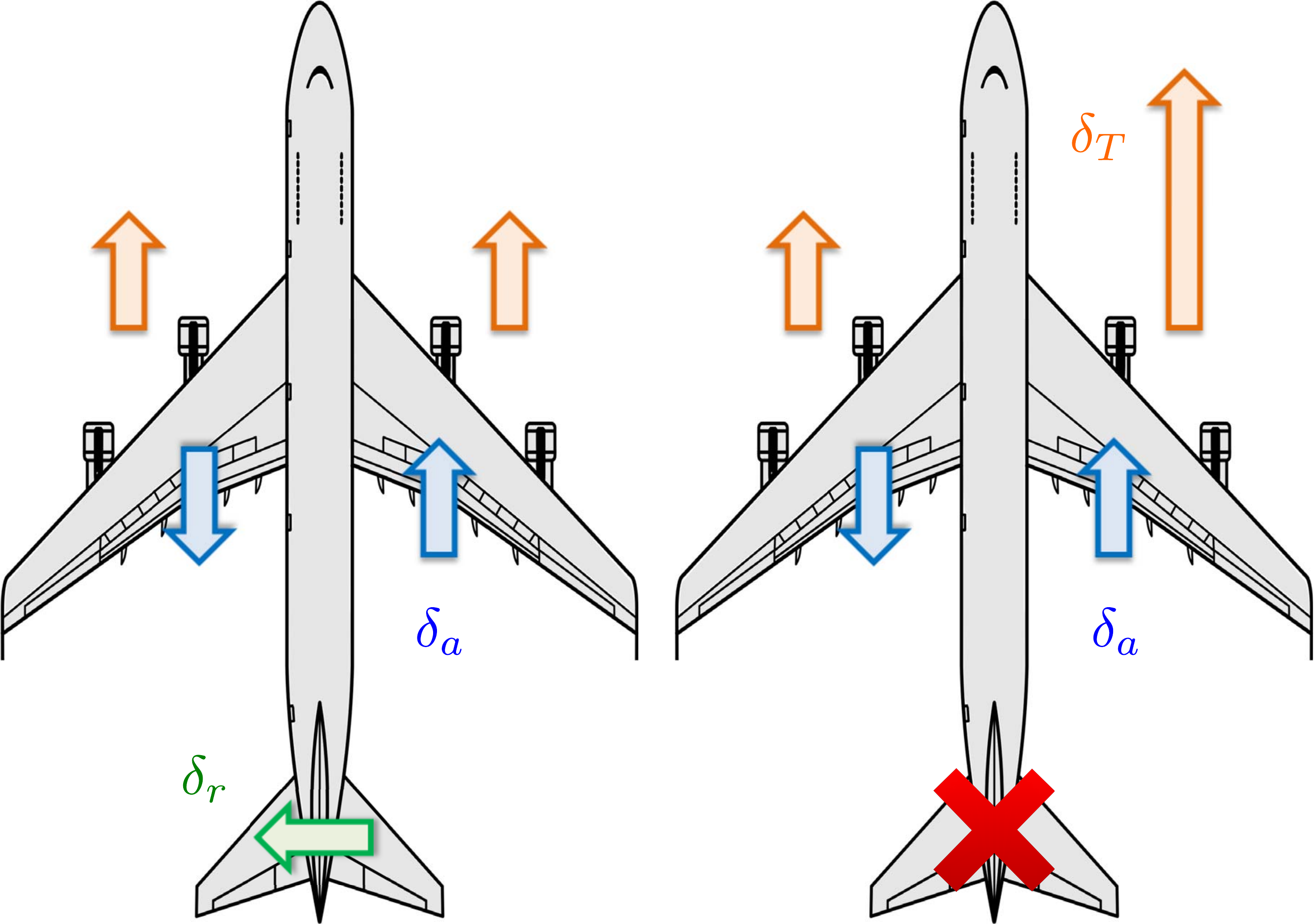}
    \caption{Top view of Boeing 747-100. Left: fault-free mode, where aileron and rudder deflections are control inputs. Right: faulty mode, where aileron deflection and differential thrust are control inputs.}
    \label{fig:Example}
\end{figure}

In order to demonstrate the effectiveness of the proposed scheme, in this section, we simulate Boeing 747-100 airplane, shown in Fig. \ref{fig:Example}. In this example, the goal is, first, to effectively detect the loss of vertical stabilizer, and then utilize a differential thrust to maintain airplane lateral/directional stability. Note that losing vertical stabilizer is a real problem in commercial airplanes. Notable examples are: 1) Japan Airlines Flight 123 in 1985, with 520 fatalities, and 2) American Airlines Flight 587 in 2001, with 265 fatalities.

Suppose that Boeing 747-100 airplane is flying at Mach 0.65 (with the corresponding airspeed of 673 [ft/sec]) at $2\times10^4$ [ft] altitude. Let $x=[\theta_r~\Delta\theta_r~\theta_s~\Delta\theta_y]^\top$ be the state of the system, where $\theta_r$ is the roll angle, $\Delta\theta_r$ is the roll rate, $\theta_s$ is the side-slip angle, and $\Delta\theta_y$ is the yaw rate. The control input is $u=[\delta_a~\delta_r~\delta_T]^\top$, where $\delta_a$ is the aileron deflection, $\delta_r$ is the rudder deflection, and $\delta_T$ is the differential thrust. The model of the system with sampling time 0.2 [s] is in the form of \eqref{eq:modesopenloop}, where $\mu_1$ and $\mu_2$ indicate, respectively, the fault-free and faulty modes. System matrices can be found in \cite{Nguyen2010} and \cite{Lu2018}. Note that the gap between the fault-free and faulty systems (computed by the \texttt{gapmetric} function in MATLAB) is 1, which means that the systems are far apart.

%with {\tiny
%\begin{eqnarray}
% A_{\mu_1}&=&\left[\begin{matrix}1 & 0.18 & -0.05 & 0.01\\-0.003 & 0.84 & -0.49 & 0.11 \\ 0.01 & 0.001 & 0.96 & -0.19\\0.001 & -0.004 & 0.2 & 0.93\end{matrix}\right],~B_{\mu_1}=\left[\begin{matrix}0.004 & 0.002 & 6.7\times10^{-4} \\0.04 & 0.02 & 0.001 \\ -2\times10^{-4} & 0.02 & -0.013\\0.002 & -0.13 & 0.0131 \end{matrix}\right],\nonumber\\
% A_{\mu_2}&=&\left[\begin{matrix}1 & 0.18 & -0.05 & 0.005\\-0.003 & 0.84 & -0.51 & 0.07 \\ 0.01 & 0.001 & 1 & -0.2\\4.2\times10^{-6} & -0.004 & 0.001 & 1\end{matrix}\right],~B_{\mu_2}=\left[\begin{matrix}0.004 & 0 & 6.7\times10^{-4} \\0.04 & 0 & 0.001 \\ -2\times10^{-4} & 0 & -0.013 \\ 0.002 & 0 & 0.0131 \end{matrix}\right],\nonumber
%\end{eqnarray}
%}
%and We use the following control laws to stabilize the system in each mode:{\tiny
%\begin{eqnarray}
% K_{\mu_1}=&&\left[\begin{matrix}-4.9688 & -6.4282 & 6.8763 & 0.2714 \\ 1.9567 & 1.2056 & 0.9075 & 1.4289 \\ 0&0&0&0\end{matrix}\right],~G_{\mu_1}=\left[\begin{matrix}4.8993 & 4.4419 \\ -2.045 & 1.0461 1.4289\\0&0\end{matrix}\right],\nonumber\\
 %K_{\mu_2}=&&\left[\begin{matrix}-3.6775 & -5.6295 & 7.1628 & 2.2242 \\0&0&0&0 \\ -1.8687 & -1.2005 & 0.824 & -2.0814 \end{matrix}\right],~G_{\mu_2}=\left[\begin{matrix}3.5497 & 5.1593 \\ 0&0 \\ 1.9685 & -1.0386 1.4289\end{matrix}\right],\nonumber
%\end{eqnarray}}

Let $H_{\omega_{\mu_1}}=H_{\omega_{\mu_2}}=2\times10^{-2}I_4$, $H_{\xi}=2\times10^{-2}I_2$, and let $r=[4.8~1.8]^\top$ be the desired reference. We consider the following constraints:
\begin{eqnarray}
&&\big\vert\mathbb{E}[\delta_a(t)]\big\vert\leq21 \text{~[deg]},~\big\vert\mathbb{E}[\delta_r(t)]\big\vert\leq3.3 \text{~[deg]},\nonumber\\
&&\big\vert\mathbb{E}[\delta_T(t)]\big\vert\leq5.2\times10^{4}\text{~[N]},\nonumber\\
&&\big\vert\mathbb{E}[\delta_T(t)-\delta_T(t-1)]\big\vert\leq2.2\times10^{4}\text{~[N]},\nonumber\\
&&\mathbb{P}(\left\vert\theta_r\right\vert\leq5^\circ)\geq0.95,~\mathbb{P}(\left\vert\theta_s\right\vert\leq2^\circ)\geq0.95,\nonumber
\end{eqnarray}
and we assume that some of these constraints can be extended for 25 steps, as follows:
\begin{eqnarray}
&&\big\vert\mathbb{E}[\delta_a(t)]\big\vert\leq25 \text{~[deg]},~\big\vert\mathbb{E}[\delta_T(t)]\big\vert\leq6.9\times10^{4}\text{~[N]},\nonumber\\
&&\mathbb{P}(\left\vert\theta_r\right\vert\leq6^\circ)\geq0.95,~\mathbb{P}(\left\vert\theta_s\right\vert\leq3^\circ)\geq0.95.\nonumber
\end{eqnarray}

\subsection{Effectiveness of the AORG}
\begin{figure}[!t]
\centering
\includegraphics[width=4cm]{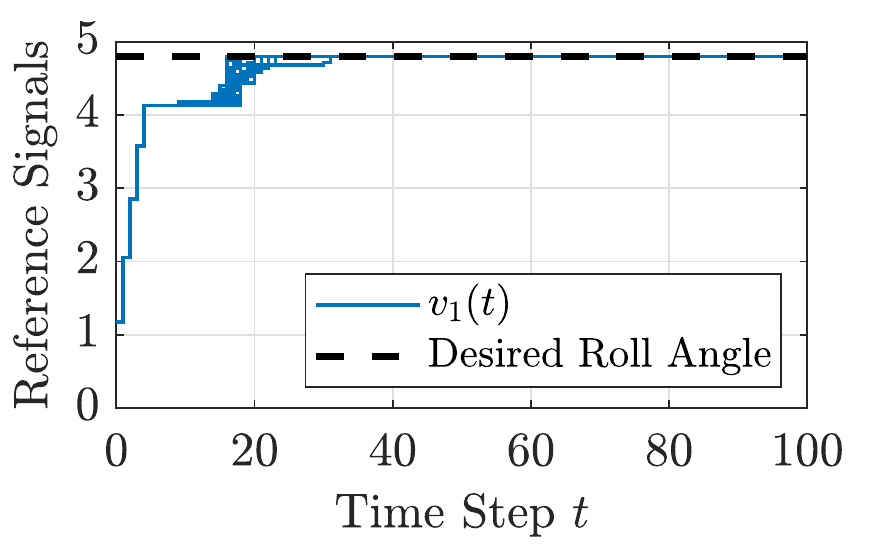}\includegraphics[width=4cm]{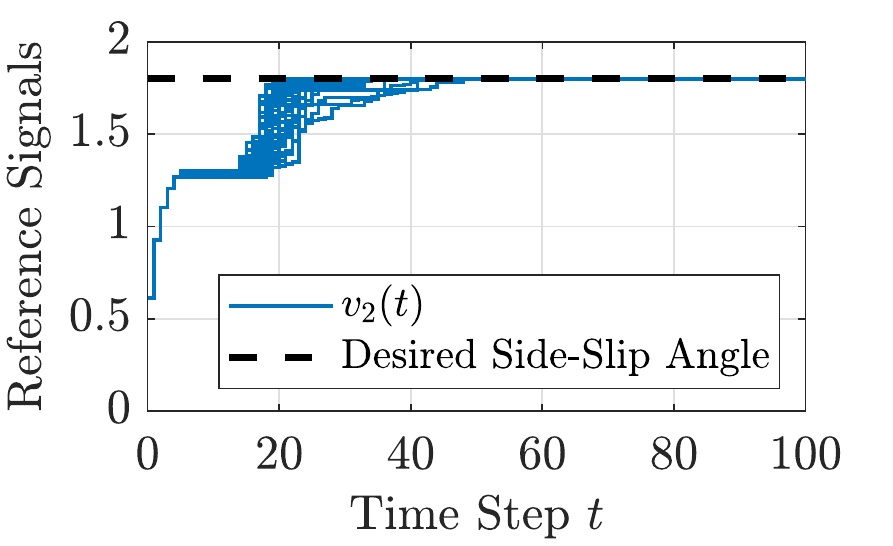}\\
\includegraphics[width=4cm]{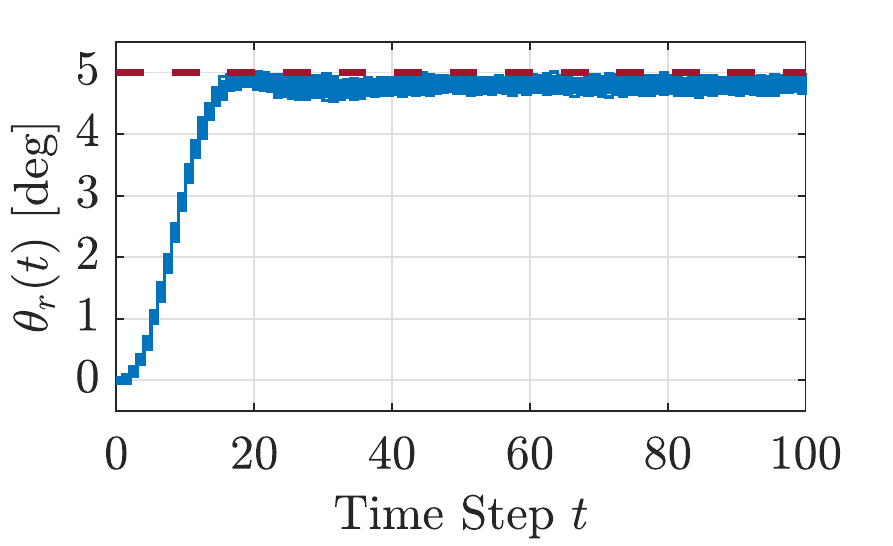}\includegraphics[width=4cm]{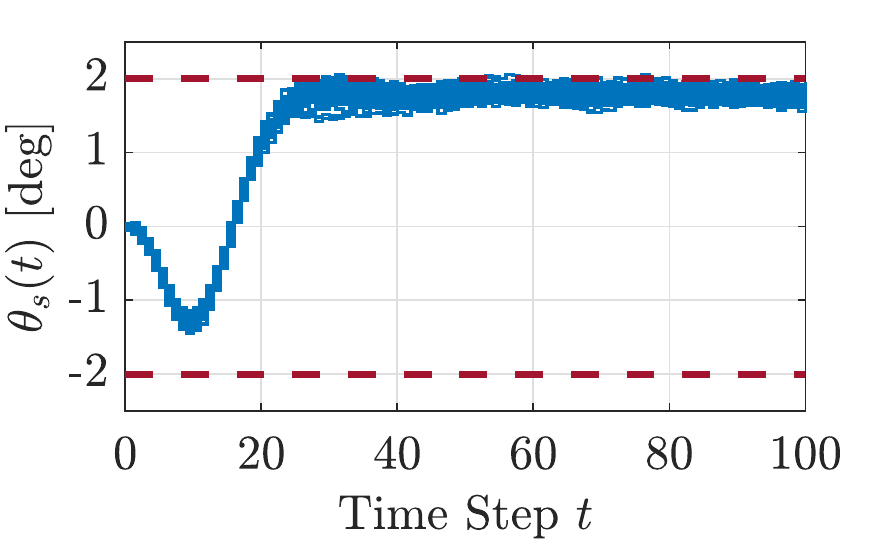}\\
\includegraphics[width=4cm]{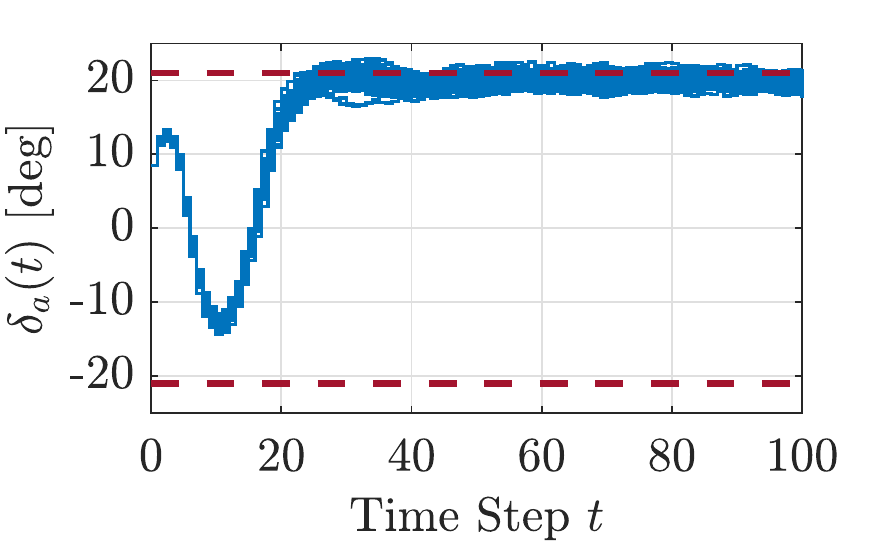}\includegraphics[width=4cm]{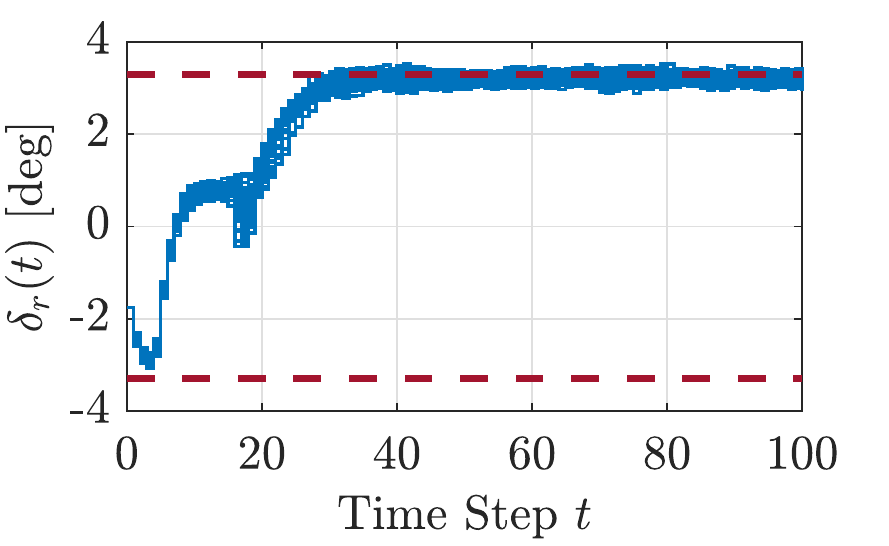}\\
\caption{Simulation results of 100 experiments conducted by the AORG, starting from the trim condition.}
\label{fig:AORG}
\end{figure}

In order to show the effectiveness of the AORG, we have run 100 simulations starting from the trim condition, i.e., $x(0)=[0~0~0~0]^\top$ and $v(0)=[0~0]^\top$, and with $T=5$. Simulation results are shown in Fig. \ref{fig:AORG}. As seen in this figure, the AORG guarantees convergence properties, while constraints are satisfied at all times.

\begin{figure}[!t]
\begin{floatrow}
\ffigbox{\includegraphics[width=4cm]{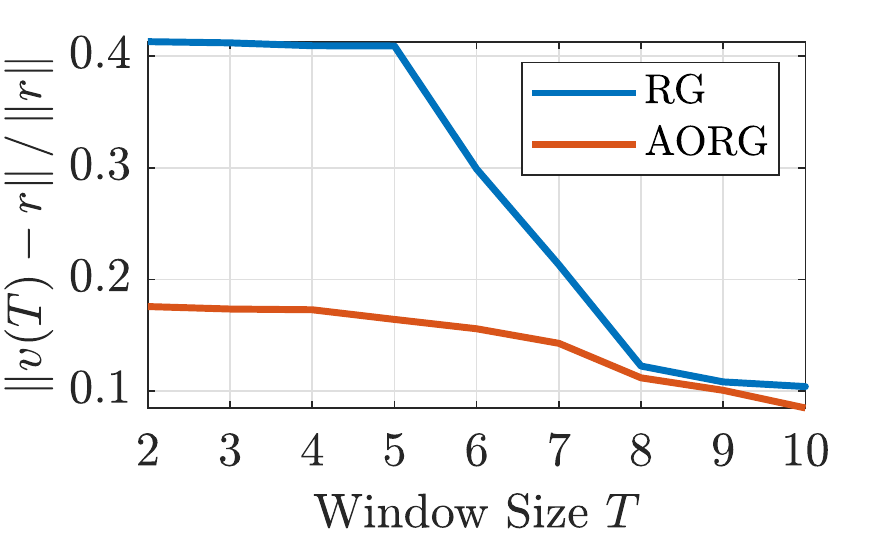}}{\caption{Comparing RG and AORG: the mean relative reference tracking error.}\label{fig:Comparison1}}
\hspace{-0.3cm}\ffigbox{\includegraphics[width=4cm]{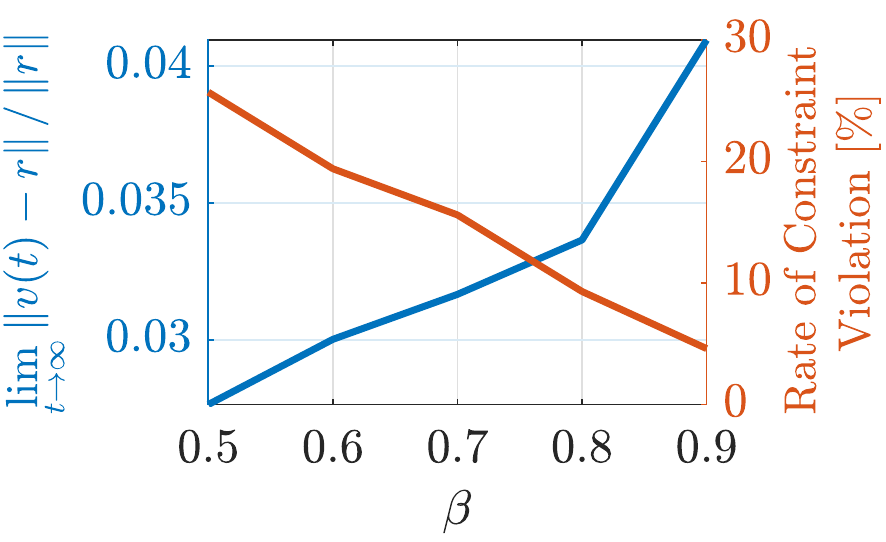}}{\caption{Impact of $\beta$ on the convergence error and the rate of constraint violation.}\label{fig:Comparison2}}
\end{floatrow}
\end{figure}

As discussed in Remark \ref{remark:AORGProperties}, 
the AORG provides a better solution compared to conventional RGs. This fact is shown in Fig. \ref{fig:Comparison1}, where the AORG is employed to compute the AR over the interval $[0,T]$. This figure reports the mean relative reference tracking error of 1000 experiments. As seen in this figure, compared to the conventional RG, the AR obtained by the AORG is closer to the desired reference at the end of the aforementioned interval, i.e., at time $T$. However, as $T$ increases, the difference between two schemes reduces. The main reason is that as time passes, the AR obtained by either schemes gets closer to the desired reference, which reduces the relative tracking error.

According to \eqref{eq:constrainttightening}, increasing $\beta$ can be interpreted as tightening the constraint. To study this, we relax the expectation constraints and we assume that $r=[5~2]^\top$. The impact of $\beta$ on the convergence error is shown in Fig. \ref{fig:Comparison2}, obtained from 1000 simulation runs. As seen in this figure, as the value $\beta$ increases, the convergence error increases as well. Also, the rate of constraint violation is around $(1-\beta)/2$. This result is expected, as the half of noises should will cause constraint violation.

\subsection{Effectiveness of the Detection Scheme}
In this section, we will assess the performance of the AORG combined with the MMAE  discussed in Section \ref{sec:controlduringtransient} in identifying the actual operating mode of the system. To do so, we have run 1000 simulations with $T_d=6$, and with random initial conditions for both fault-free and faulty modes. Results are reported in Table \ref{tab:comparison}. As seen in this table, when the detection objective function is not taken into account in determining the AR sequence, the MMAE  can identify the actual operating mode of the system with 98.7\% of correctness when the vertical stabilizer works faultlessly. However, when the vertical stabilizer fails, the actual operating mode of the system can be identified with 58.1\% of correctness.

As discussed in Section \ref{sec:controlduringtransient}, considering the detection objective function in determining the AR sequence can improve the detection performance. This fact can be seen in Table \ref{tab:comparison}. As reported in this table, by simultaneously considering the detection and control objective functions, the MMAE  can identify the actual operating mode with 97.6\% of correctness for the fault-free mode, and with 100\% for the faulty mode. Note that the value of the control objective function (i.e., $\sum_i\sum_j\kappa_j^i$) with considering the detection objective function is 34.4\% less than that of the case where the detection objective function is not considered.

{\color{black}For comparison purposes, we implement the fault-tolerant Model Predictive Controller (MPC) described in \cite{Camacho2010}. As seen in Table \ref{tab:comparison}, when we use the MPC scheme, the MMAE  can identify the actual operating mode of the system with 83.9\% of correctness when the vertical stabilizer works faultlessly, and with 72.4\% when the vertical stabilizer fails.}

\begin{table}[!t]
\centering
\caption{Effectiveness of the AORG scheme combined with the MMAE in detecting faults.} \label{tab:comparison}
\begin{tabular}{c|c|c}
\hline
Controller & Mode & Correct Identification\\
& & (Percentage) \\
\hline
{\color{black}MPC \cite{Camacho2010}} & {\color{black}Fault-Free} & {\color{black}83.9\%} \\
& {\color{black}Faulty} & {\color{black}72.4\%} \\
\hline
AORG & Fault-Free & 58.1\% \\
& Faulty & 97.6\%\\
\hline
AORG  & Fault-Free & 98.7\% \\
with Detection & Faulty & 100\%\\
\hline
\end{tabular}
\end{table}

\begin{figure}[!t]
\centering
\includegraphics[width=\columnwidth]{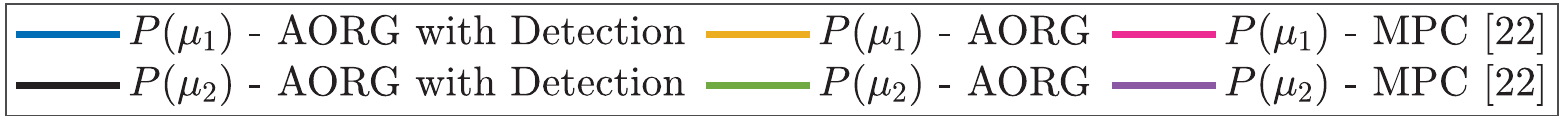}\\
\includegraphics[width=4cm]{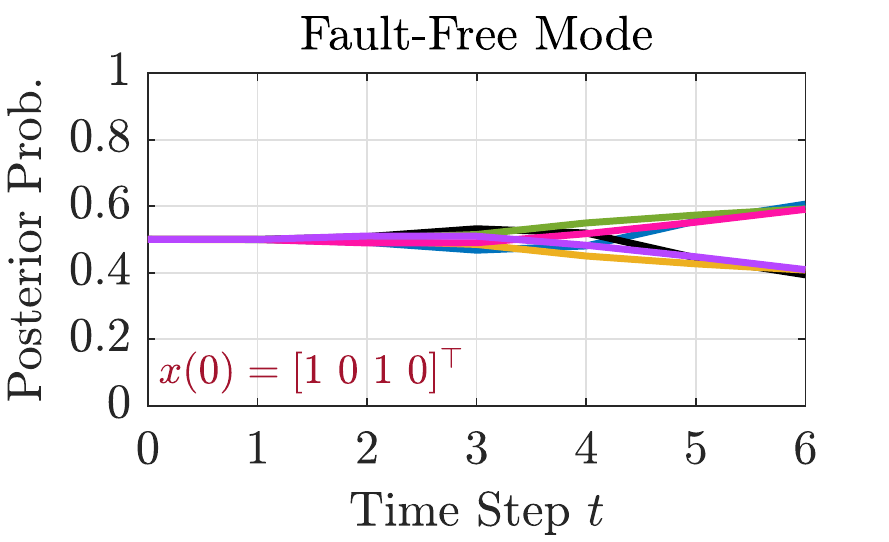}\includegraphics[width=4cm]{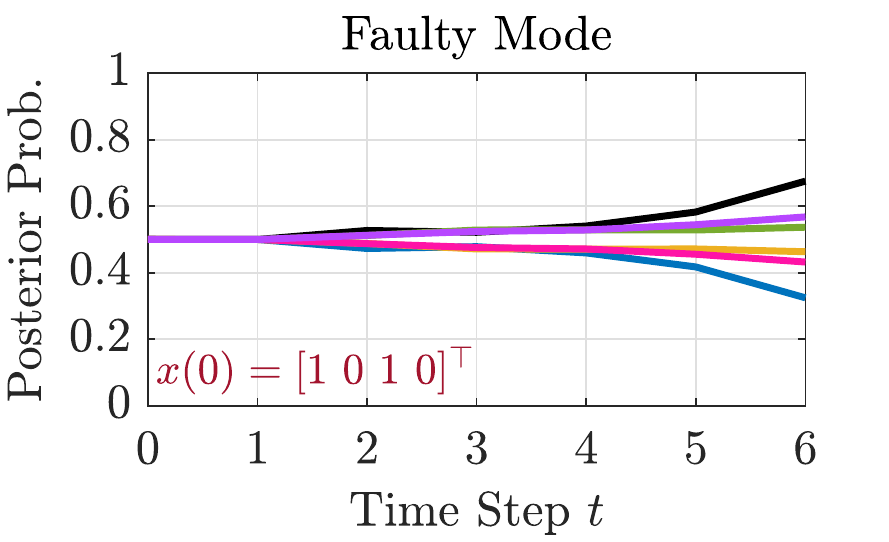}\\
\includegraphics[width=4cm]{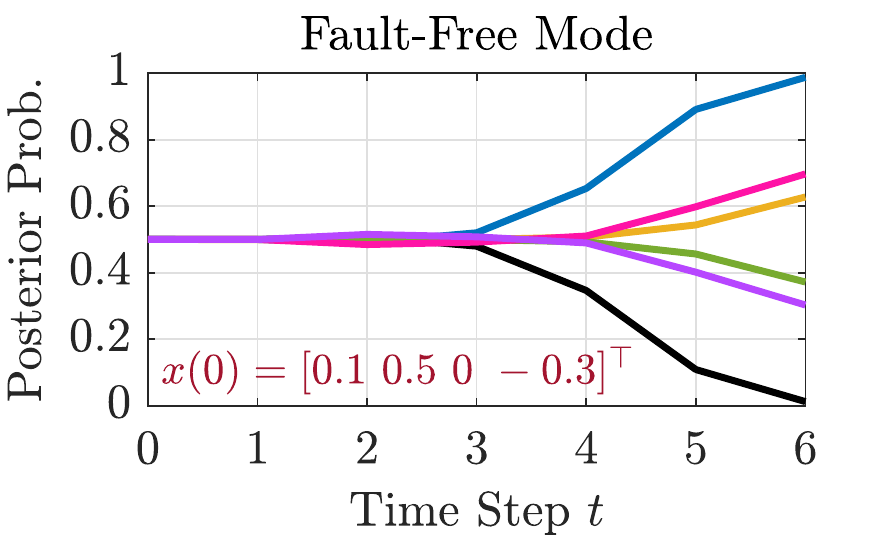}\includegraphics[width=4cm]{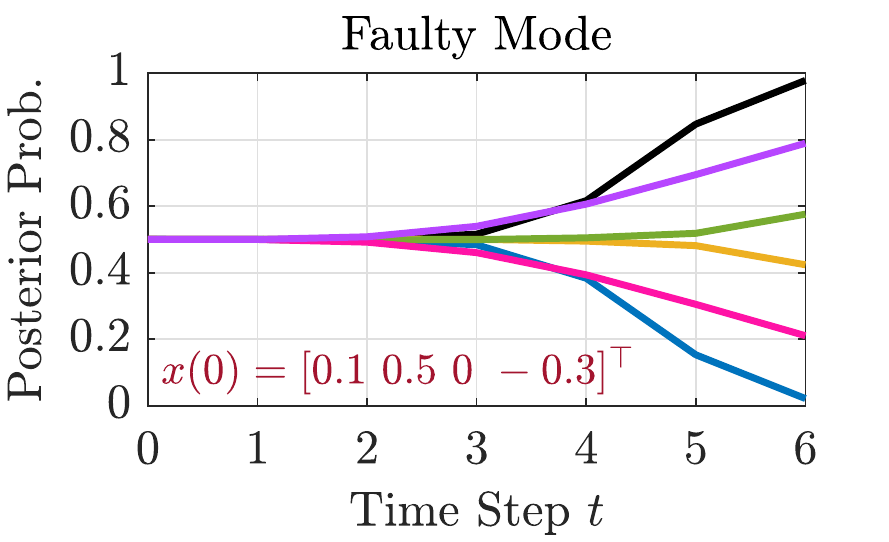}
\caption{{\color{black}Posterior probabilities obtained by the MMAE  for two typical initial conditions.}}
\label{fig:Detection}
\end{figure}

The posterior probabilities obtained by the MMAE  for two typical initial conditions are shown in Fig. \ref{fig:Detection}. As seen in Fig. \ref{fig:Detection}, when the initial condition is $x(0)=[1~0~1~0]^\top$ and when the system operates under the fault-free mode, without considering the detection objective function in determining the AR sequence, the MMAE  does not identify the actual operating mode of the system. However, by taking into account the detection objective function, the MMAE  identifies the actual operating mode of the system, though with a low level of confidence.

\subsection{Effectiveness of the Reconfiguration Scheme}
\begin{figure}[!t]
\centering
\includegraphics[width=4cm]{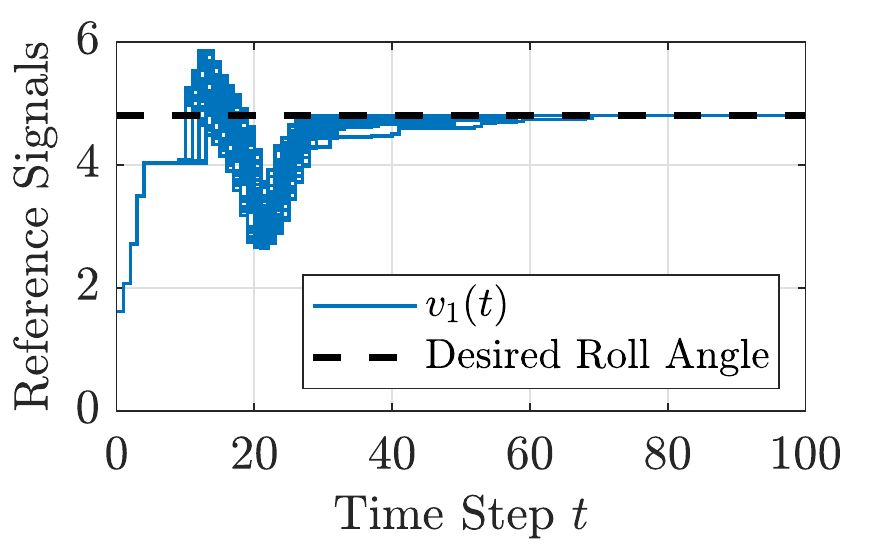}\includegraphics[width=4cm]{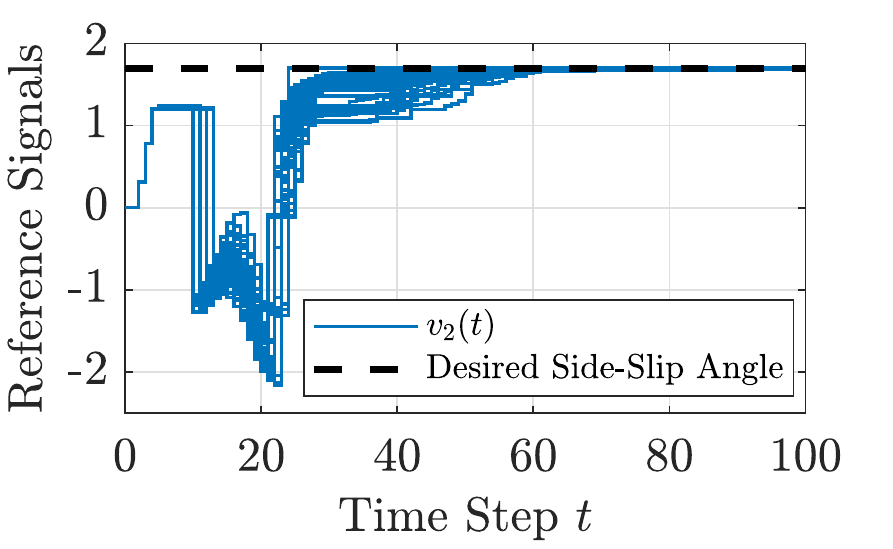}\\
\includegraphics[width=4cm]{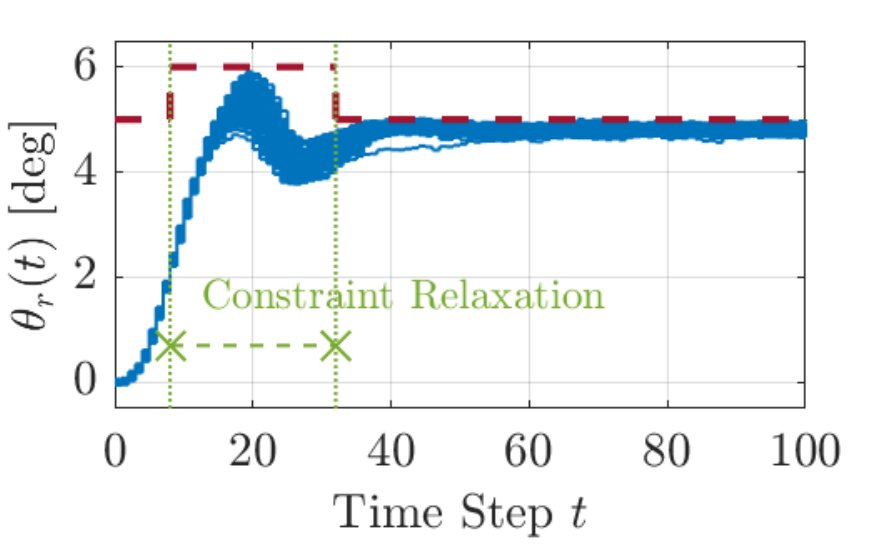}\includegraphics[width=4cm]{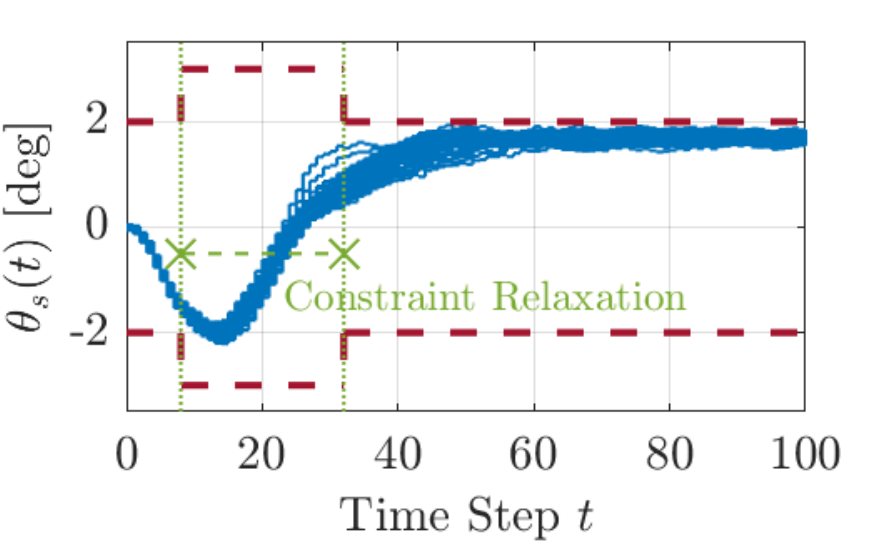}\\
\includegraphics[width=4cm]{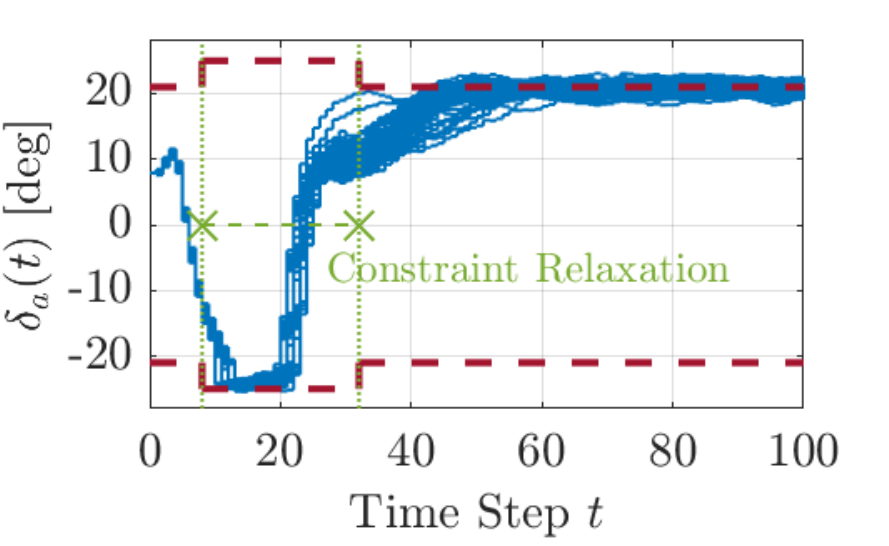}\includegraphics[width=4cm]{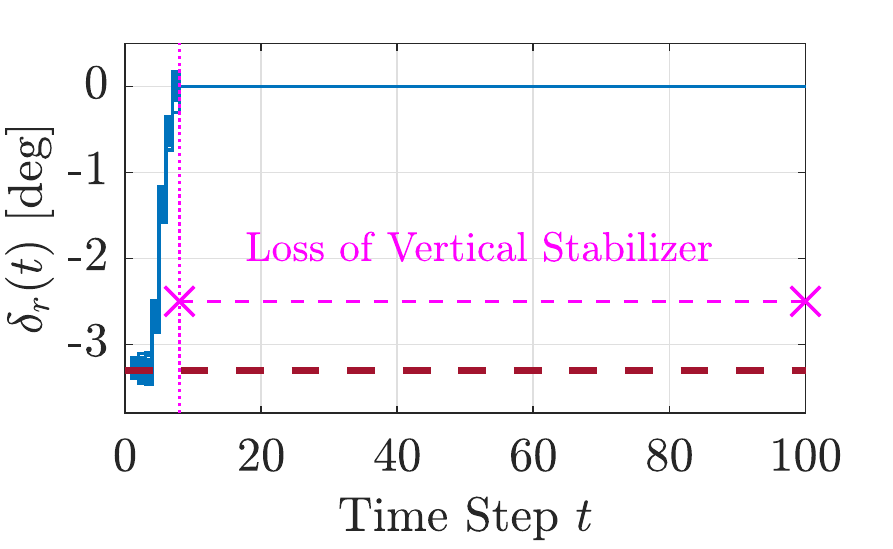}\\\includegraphics[width=4cm]{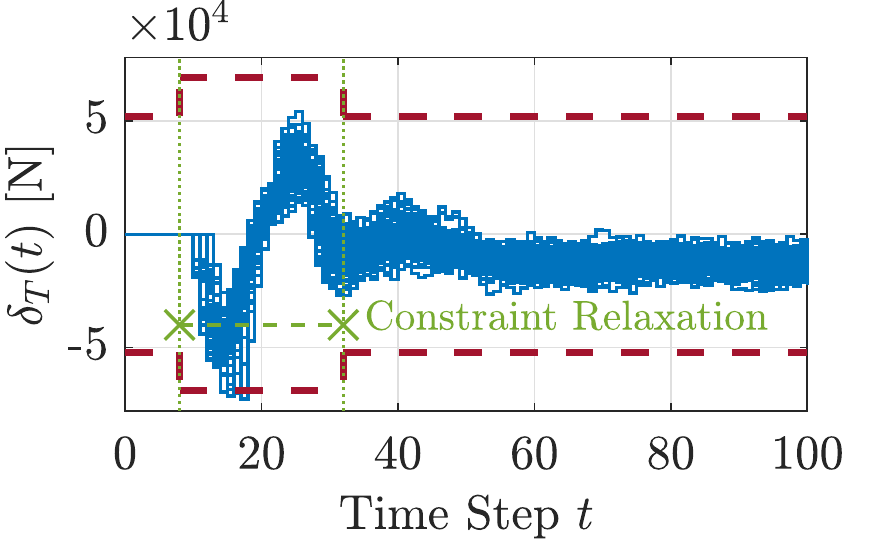}\\
\caption{Simulation results of 100 experiments started from the trim condition. The vertical stabilizer fails at $t=8$, and the reconfiguration scheme employs the differential thrust to recover tracking and constraint satisfaction properties.}
\label{fig:Recovery}
\end{figure}

In this section, we will evaluate the effectiveness of the proposed reconfiguration scheme. Suppose that the vertical stabilizer fails at $t=8$. Suppose that $T_d=6$ and $T_r=13$. Simulation results for 100 experiments starting from the trim condition are shown in Fig. \ref{fig:Recovery}. As seen in this figure, once the vertical stabilizer fails and the MMAE detects this failure, the reconfiguration schemes gets involved by computing a AR sequence via Lemma \ref{lemma:reconfiguration}. As shown in Fig. \ref{fig:Recovery}, the computed AR sequence can safely recover the tracking and constraint satisfaction properties.

%%%%%%%%%%%%%%%%%%%
\section{Conclusion}\label{sec:conclusion}
This paper proposed a fault-tolerant constrained control scheme. First, a new RG-based constrained control scheme, called AORG, was presented. The main feature of this scheme is that it computes the AR sequence for an interval entirely at the beginning of the interval. It's convergence and constraint-handling properties are proven rigorously. The AORG was combined with the MMAE which was used to detect the fault. It was shown that the AR sequence can be determined such that the detection performance is optimized, while enforcing constraints satisfaction at all times and ensuring reference tracking. Finally, a reconfiguration scheme based on recoverable sets was presented, which can maintain system viability and functionality despite the presence of the fault. The effectiveness of the proposed scheme is validated through extensive simulation studies carried out on Boeing 747-100. {\color{black}As future work, we plan to extend our method to deal with non-Gaussian and possibly correlated noise.}

%---------------------------------------------
%---------------------------------------------
%\balance
\bibliographystyle{ieeetr}
\bibliography{ref}
%%%%%%%%%%%%%%%%%%%%%%%%%%%%%%%%%%%%%%%%%%%%%%%%%%%%%%%%%%%%%%%%%%%%%%%%%%%%%%%%
\parpic{\includegraphics[width=1in,height=1in,clip,keepaspectratio]{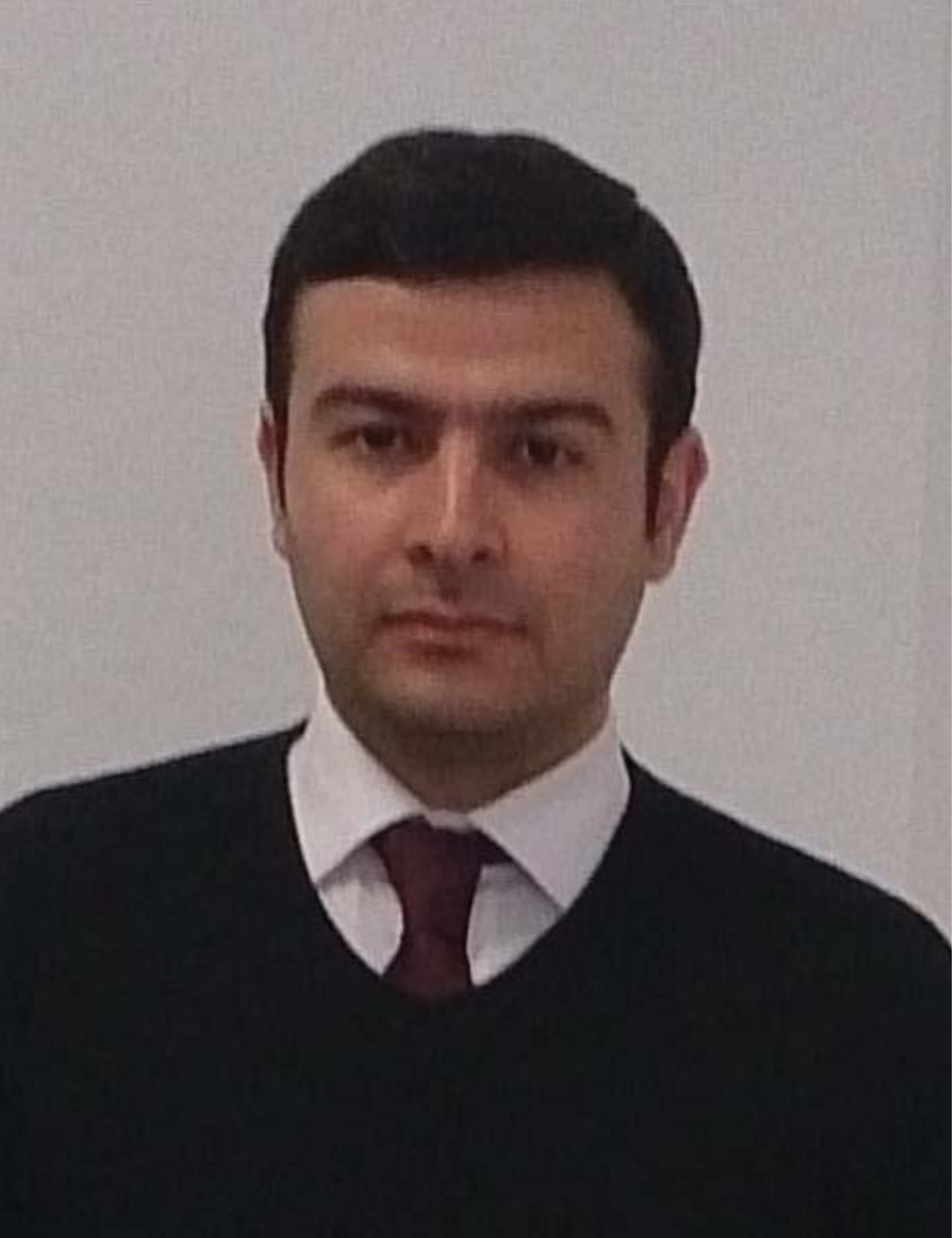}}
\textbf{Mehdi Hosseinzadeh} received his Ph.D. degree in Electrical Engineering-Control from the University of Tehran, Iran, in 2016. From 2017 to 2019, he was a postdoctoral researcher at Universit\'{e} Libre de Bruxelles, Brussels, Belgium. In 2018, he was a visiting researcher at University of British Columbia, Canada. He is currently a postdoctoral research associate at Washington University in St. Louis, MO, USA. His research interests include nonlinear and adaptive control, constrained control, and safe and robust control of autonomous systems.

\parpic{\includegraphics[width=1in,height=1in,clip,keepaspectratio]{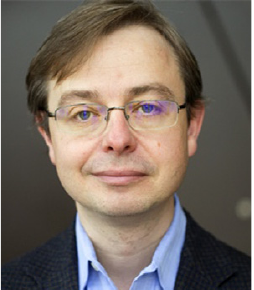}}
\textbf{Ilya Kolmanovsky}  is a professor in the department of aerospace engineering at the University of Michigan,with research interests in control theory for systems with state and control constraints, and in control applications to aerospace and automotive systems. He received his Ph.D. degree in aerospace engineering from the University of Michigan in 1995. Prof. Kolmanovsky is a Fellow of IEEE and is named as an inventor on 104 United States patents.

\parpic{\includegraphics[width=1in,height=1in,clip,keepaspectratio]{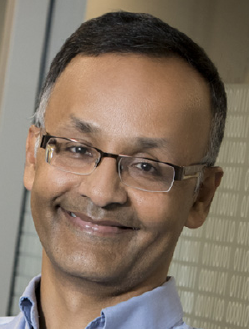}}
\textbf{Sanjoy Baruah} joined Washington University in St. Louis in September 2017. He was previously at the University of North Carolina at Chapel Hill (1999\textendash 2017) and the University of Vermont (1993\textendash 1999). His research interests and activities are in real-time and safety-critical system design, scheduling theory, resource allocation and sharing in distributed computing environments, and algorithm design and analysis.

\parpic{\includegraphics[width=1in,height=1in,clip,keepaspectratio]{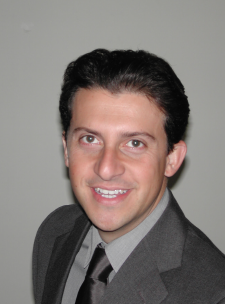}}

\textbf{Bruno  Sinopoli} received his Ph.D. in Electrical Engineering from the University of  California at Berkeley, in 2005. After a postdoctoral position at Stanford University, he was the faculty at Carnegie Mellon University from 2007 to 2019, where he was full  professor  in  the  Department  of  Electrical  and Computer  Engineering  with  courtesy  appointments in  Mechanical  Engineering  and  in  the  Robotics  Institute  and  co-director  of  the  Smart  Infrastructure Institute. In  2019  he  joined Washington University in Saint Louis, where he is the chair of the Electrical and Systems Engineering department. He was awarded the 2006 Eli Jury  Award  for  outstanding  research  achievement  in  the  areas  of  systems, communications,  control  and  signal  processing  at  U.C.  Berkeley,  the  2010 George Tallman Ladd Research Award from Carnegie Mellon University and the NSF Career award in 2010. His research interests include the modeling,analysis  and  design  of  Secure  by  Design  Cyber-Physical  Systems  with applications  to  Energy  Systems,  Interdependent  Infrastructures  and  Internet of Things.

\end{document}